\documentclass[a4paper]{article}
\usepackage[dvipdfmx]{graphicx}
\usepackage[top=2cm, bottom=2cm, left=2cm, right=2cm]{geometry}
\usepackage{amsthm}
\usepackage{newtxtext, newtxmath}
\usepackage{cite}
\usepackage[linesnumbered,lined,ruled,commentsnumbered]{algorithm2e}
\usepackage{float}
\usepackage{color} 
\usepackage{comment} 

\usepackage{hyperref} 
\usepackage{xcolor}
\hypersetup{
    colorlinks=true,
    citecolor=blue,
    linkcolor=blue,
    urlcolor=orange,
}

\theoremstyle{definition}
\newtheorem{dfn}{Definition}[section]
\newtheorem{thm}[dfn]{Theorem}
\newtheorem{conj}[dfn]{Conjecture}
\newtheorem{lem}[dfn]{Lemma}

\let\oldenumerate\enumerate
\renewcommand{\enumerate}{
   \oldenumerate
   \setlength{\itemsep}{0cm}
   \setlength{\parskip}{0cm}
}
\let\olditemize\itemize
\renewcommand{\itemize}{
   \olditemize
   \setlength{\itemsep}{0cm}
   \setlength{\parskip}{0cm}
}

\title{Improved 2-distance coloring of planar graphs with maximum degree 5}
	\author{Kengo Aoki\thanks{Department of Informatics, Graduate School of Integrated Science and Technology, Shizuoka University, 3-5-1 Johoku, Naka-ku, Hamamatsu-shi, 432-8011, Shizuoka, Japan, aoki.kengo.19@shizuoka.ac.jp}}
\date{}

\begin{document}
\maketitle

\begin{abstract}
A 2-distance $k$-coloring of a graph $G$ is a proper $k$-coloring such that any two vertices at distance two or less get different colors.
The 2-distance chromatic number of $G$ is the minimum $k$ such that $G$ has a 2-distance $k$-coloring, denote as $\chi_2(G)$.
In this paper, we show that $\chi_2(G) \leq 17$ for every planar graph $G$ with maximum degree $\Delta \leq 5$, which improves a former bound $\chi_2(G) \leq 18$.
\end{abstract}

\section{Introduction}
We use notations based on \cite{MR2368647}, \cite{chen20222} and \cite{zhu2021wegner}.
In this paper, all graphs are simple, finite, and planar.
For a graph $G$, we denote the set of vertices, the set of edges and the set of faces by $V(G)$, $E(G)$ and $F(G)$, respectively.
The ends of an edge are said to be \textit{incident} with edge, and \textit{vice versa}. 
Two vertices which are incident with a common edge are \textit{adjacent}, and two distinct adjacent vertices are \textit{neighbours}.
The set of neighbours of a vertex $v$ in a graph $G$ is denoted by $N_G(v)$.
If $e$ is an edge of $G$, we may obtain a graph on $|E(G)| - 1$ edges by deleting $e$ from $G$ but leaving the vertices and remaining edges intact.
The resulting graph is denoted by $G - e$.
Similarly, if $v$ is a vertex of $G$, we may obtain a graph on $|V(G)| - 1$ vertices by deleting from $G$ the vertex $v$ together with all the edges incident with $v$.
The resulting graph is denoted by $G - \{v\}$.
The \textit{degree} of a vertex $v$ in a graph $G$, denoted by $\mathrm{deg}_G(v)$, is the number of edges of $G$ incident with $v$.
The maximum degree and minimum degree of a graph $G$ are respectively denoted by $\Delta(G)$ and $\delta(G)$.
A vertex of degree $k$ (respectively, at least $k$, at most $k$) is said to be $k$-\textit{vertex}(respectively, $k^+$-vertex, $k^-$-vertex).
A face is said to be \textit{incident} with the vertices and edges in its boundary, and two faces are \textit{adjacent} if their boundaries have an edge in common.
The \textit{degree} of a face $f$ in a graph $G$, denoted by $\mathrm{deg}_G(f)$, is the number of edges incident to $f$.
A face of degree $k$ (respectively, at least $k$, at most $k$) is said to be $k$-\textit{face}(respectively, $k^+$-face, $k^-$-face).
Let $t(v)$ be the number of 3-faces incident to a vertex $v$.
A $[v_1v_2\cdots v_k]$ is a $k$-face with vertices $v_1,v_2, \cdots ,v_k$ on its boundary.
A $(x_1,x_2,\cdots,x_k)$-\textit{face} is a $k$-face with vertices of degrees $x_1,x_2,\ldots, x_k$.
Let $\phi$ be a partial coloring of a graph $G$.
For a vertex $v$ in a graph $G$, let $C_\phi(v)$ denotes the set of colors which are assigned on the vertices within distance two to $v$.
A 2-distance $k$-coloring of a graph $G$ is a mapping $\phi :V(G) \rightarrow \{1,2,\cdots,k\}$ such that $\phi(v_1) \neq \phi(v_2)$ if any two vertices $v_1,v_2$ with $\mathrm{dist}(v_1, v_2) \leq 2$ where $\mathrm{dist}(v_1, v_2)$ is the distance between the two vertices $v_1$ and $v_2$.
The 2-distance chromatic number of G is the minimum $k$ such that G has a 2-distance $k$-coloring, denote as $\chi_2(G)$.

In 1977, Wegner make the following conjecture.
\begin{conj}\cite{wegner1977graphs}\label{wegner}
If $G$ is a planar graph, then $\chi_2(G) \leq \Delta(G)+5$ if $4\leq \Delta(G) \leq 7$ and $\chi_2(G) \leq \lfloor \frac{3\Delta(G)}{2} \rfloor + 1$ if $\Delta(G) \geq 8$.
\end{conj}
Conjecture~\ref{wegner} is still open. 
The upper bound of $\chi_2(G)$ is getting smaller with $\Delta(G) \leq 5$.
Heuvel and McGuinness \cite{van2003coloring} proved that $\chi_2(G) \leq 9\Delta(G) - 19$ if $\Delta(G) \geq 5$.
Zhu and Bu\cite{zhu2018minimum} proved that $\chi_2(G) \leq 20$.
Chen, Miao and Zhou\cite{chen20222} proved that $\chi_2(G) \leq 19$.
J.Zhu, Bu and H.Zhu\cite{zhu2021wegner} proved that $\chi_2(G) \leq 18$.
In this paper, we improve a result in \cite{zhu2021wegner} by solving the following theorem.
\begin{thm}\label{main}
If $G$ is a planar graph with maximum degree $\Delta(G) \leq 5$, then $\chi_2(G) \leq 17$.
\end{thm} 

\section{Reducible configurations}
Let $G$ be a minimum counterexample with minimum $|V(G)| + |E(G)|$ to Theorem~\ref{main}.
That is $\chi_2(G) > 17$.
The minimum means that for any subgraph $G'$ obtained from $G$ with $\Delta(G') \leq 5$ and $|V(G')| + |E(G')| < |V(G)| + |E(G)|$, the inequality $\chi_2(G') \leq 17$ holds.
Let $C = \{1,2, \cdots , 17\}$ be a set of colors.
In this section, we indicate reducible configurations.
Lemmas cited from \cite{chen20222}, \cite{zhu2018minimum} and\cite{zhu2021wegner}, respectively assumes $\chi_2(G') \leq 19$, $\chi_2(G') \leq 20$ and $\chi_2(G') \leq 18$, they also hold true even when assuming $\chi_2(G') \leq 17$.

\begin{lem}[{\cite[Lemma~2.1]{zhu2018minimum}}]\label{lem1}
A graph $G$ is connected.
\end{lem}

\begin{lem}[{\cite[Lemma~2.1]{chen20222}}]\label{lem2}
In a graph $G$, there is no cut edge.
\end{lem}

\begin{lem}[{\cite[Lemma~2.2]{chen20222}}]\label{lem3}
In a graph $G$, $\delta(G) \geq 3$.
\end{lem}
From this point onward, we use Lemma~\ref{lem1} to Lemma~\ref{lem3} without explicit citation.
\begin{lem}[{\cite[Lemma~2.3]{chen20222}}]\label{lem4}
Every 3-vertex is adjacent to three 5-vertices.
\end{lem}

\begin{lem}[{\cite[Lemma~2.4]{chen20222}}]\label{lem5}
There is no 3-vertex on 3-face in a graph $G$.
\end{lem}

\begin{figure}[ht]
\begin{center}
\includegraphics[scale=0.25]{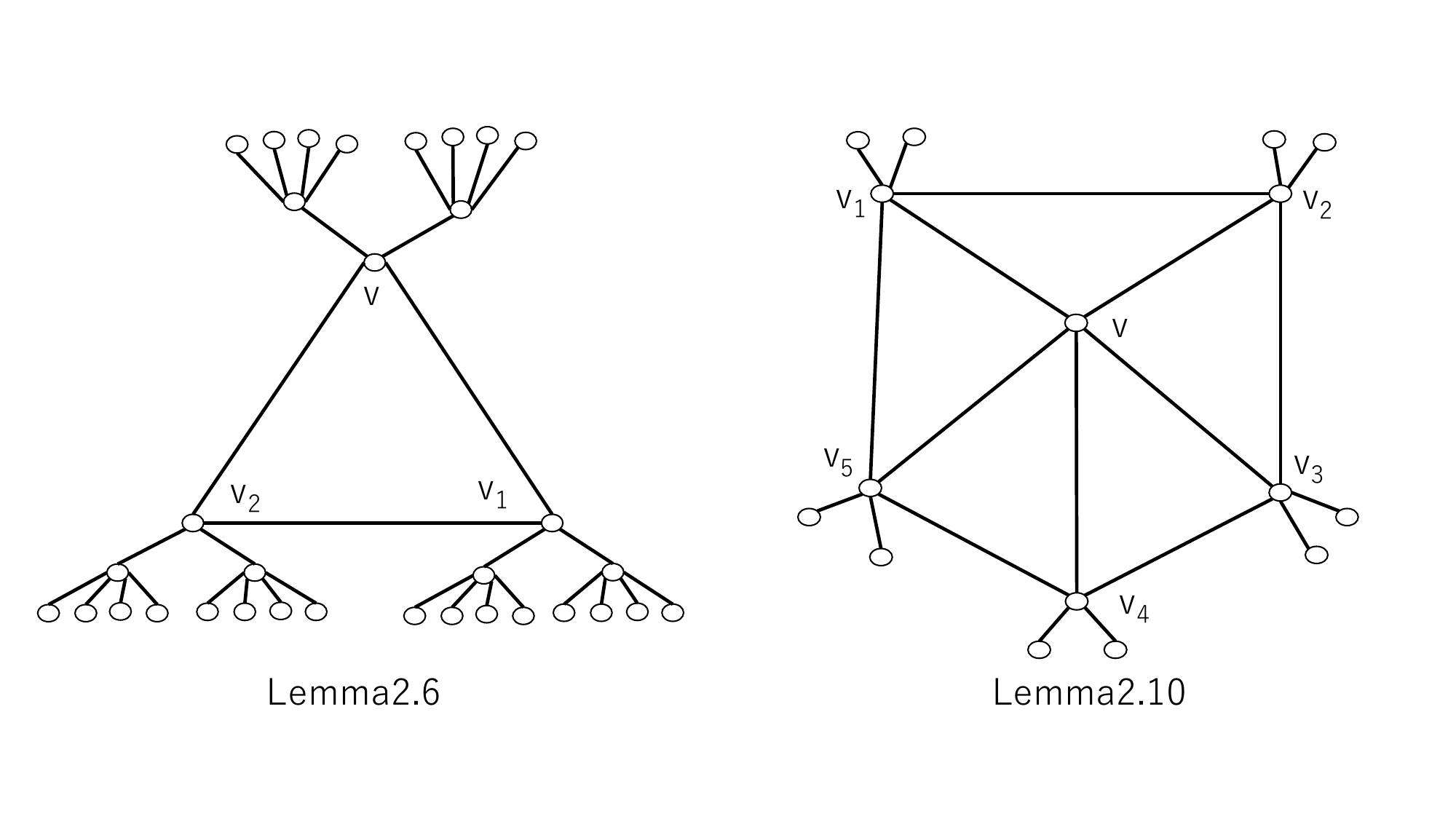}
\caption{Illustrations of Lemma~\ref{lem6} and Lemma~\ref{lem10}.}
\label{fig_lem_6_lem10}
\end{center}
\end{figure}

\begin{lem}\label{lem6}
There are no (4,4,4)-faces in a graph $G$.(See Figure~\ref{fig_lem_6_lem10}.)
\begin{proof}
Assume there exists a 3-face $[vv_1v_2]$ such that $\mathrm{deg}_G(v) =\mathrm{deg}_G(v_1) = \mathrm{deg}_G(v_2) = 4$.
Let $G' = G - vv_1$.
By the minimality of $G$, $G'$ has a 2-distance 17-coloring $\phi'$.
Every vertex in $V(G)$ is colored using $\phi'$.
Erase the color of $v$ and $v_1$, and recount the available colors for the two vertices.
Since $\Delta(G) \leq 5$, it follows that $|C_{\phi'}(v)| \leq 5 + 5 + 3 + 2 = 15$ and $|C| - |C_{\phi'}(v)| \geq 2$, $|C_{\phi'}(v_1)| \leq 5 + 5 + 3 + 2 = 15$ and $|C| - |C_{\phi'}(v_1)| \geq 2$.
If $v$ and $v_1$ are colored with $\phi'(v) \in C \setminus C_{\phi'}(v)$, $\phi'(v_1) \in C \setminus C_{\phi'}(v_1)$ and $\phi'(v) \neq \phi'(v_1)$, then $\phi'$ can be extended to a 2-distance 17-coloring of $G$, which is a contradiction.
\end{proof}
\end{lem}

\begin{lem}[{\cite[Lemma~2.10]{chen20222}}]\label{lem7}
There is at most one 3-vertex on 5-face in a graph $G$.
\end{lem}

\begin{lem}[{\cite[Lemma~3.6]{zhu2021wegner}}]\label{lem8}
Every 4-vertex is incident to at most one 3-face.
\end{lem}

\begin{lem}[{\cite[Lemma~2.6]{chen20222}}]\label{lem9}
In a graph $G$, if a 3-vertex $v$ is incident to 4-face, then other three vertices on the 4-face are 5-vertices.
\end{lem}

\begin{lem}\label{lem10}
If $v$ is a 5-vertex, then $t(v) \leq 4$.(See Figure~\ref{fig_lem_6_lem10}.)
\begin{proof}
Let $N_G(v) = \{v_1,v_2,v_3,v_4,v_5\}$.
Suppose $v$ is incident to five 3-faces $[vv_1v_2]$, $[vv_2v_3]$, $[vv_3v_4]$, $[vv_4v_5]$ and $[vv_5v_1]$.
We make graph $G'$ so that the vertices with distance less than or equal to 2 in graph $G$ also have distance less than or equal to 2 in graph $G'$.
This way is the same in the subsequent proofs of Lemmas when making $G'$ by deleting vertices from $G$.
Let $G' = G - \{v\}$.
By the minimality of $G$, $G'$ has a 2-distance 17-coloring $\phi'$.
Let $\phi$ be a coloring of $G$ such that every vertex in $V(G)$, except for $v$, is colored using $\phi'$.
Since $\Delta(G) \leq 5$, it follows that $|C_\phi(v)| \leq 3 + 3 + 3 + 3 + 3 = 15$ and $|C| - |C_\phi(v)| \geq 2$.
If $v$ is colored with $\phi(v) \in C \setminus C_\phi(v)$, then there exists a coloring $\phi$ of $G$ such that $\chi_2(G) \leq 17$, which is a contradiction.
\end{proof}
\end{lem}

\begin{figure}[ht]
\begin{center}
\includegraphics[scale=0.25]{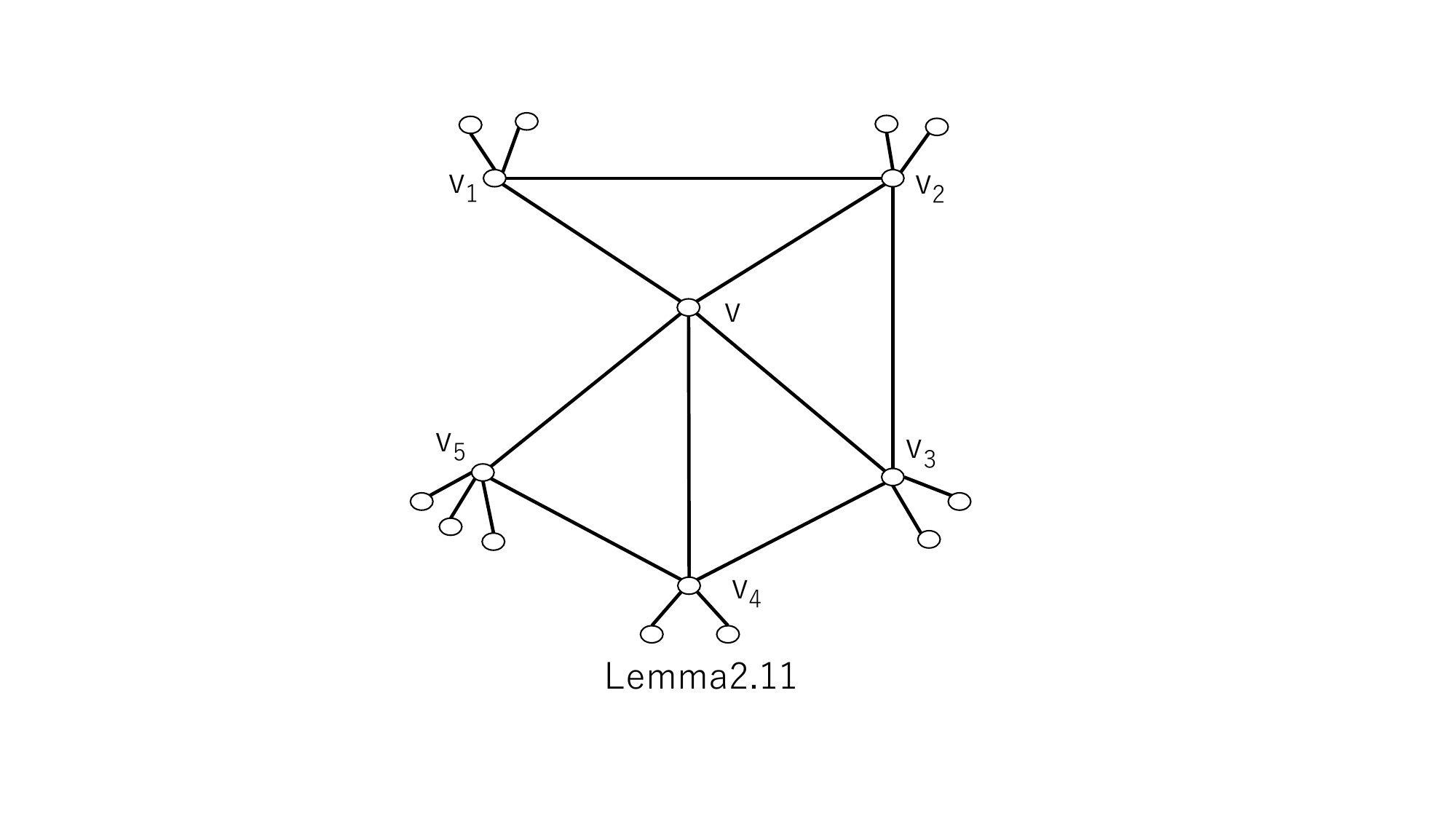}
\caption{Illustration of Lemma~\ref{lem11}: Assuming a 5-vertex $v$ is incident to three (5,5,5)-faces and one (5,5,4)-face.}
\label{fig_lem11}
\end{center}
\end{figure}

\begin{lem}\label{lem11}
If $v$ is a 5-vertex with $t(v) = 4$, then the degree of all vertices adjacent to $v$ is five.(See Figure~\ref{fig_lem11}.)
\begin{proof}
Let $N_G(v) = \{v_1,v_2,v_3,v_4,v_5\}$.
Assume $v$ is a 5-vertex which is incident to four 3-faces $[vv_1v_2]$, $[vv_2v_3]$, $[vv_3v_4]$ and $[vv_4v_5]$.
By Lemma~\ref{lem5} and Lemma~\ref{lem8}, $\mathrm{deg}_G(v_2) =\mathrm{deg}_G(v_3) = \mathrm{deg}_G(v_4) = 5$.
Suppose $\mathrm{deg}_G(v_1) = 4$.
Let $G' = G - \{v\} + v_1v_5$.
By the minimality of $G$, $G'$ has a 2-distance 17-coloring $\phi'$.
Let $\phi$ be a coloring of $G$ such that every vertex in $V(G)$, except for $v$, is colored using $\phi'$.
Since $\Delta(G) \leq 5$, it follows that $|C_\phi(v)| \leq 3 + 3 + 3 + 3 + 4 = 16$ and $|C| - |C_\phi(v)| \geq 1$.
If $v$ is colored with $\phi(v) \in C \setminus C_\phi(v)$, then there exists a coloring $\phi$ of $G$ such that $\chi_2(G) \leq 17$, which is a contradiction.
\end{proof}
\end{lem}

\begin{figure}[ht]
\begin{center}
\includegraphics[scale=0.25]{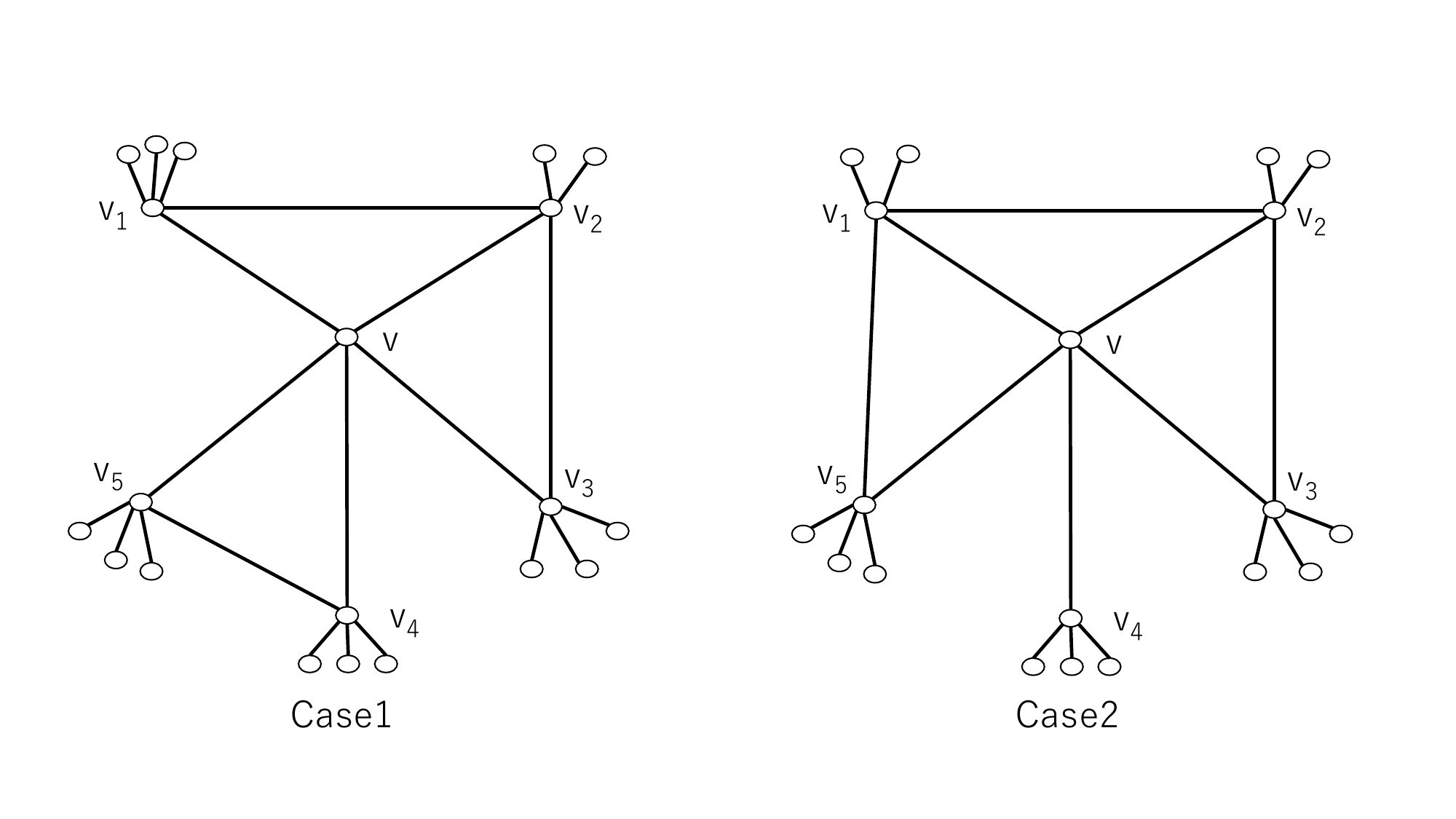}
\caption{Illustrations of Lemma~\ref{lem12}: There are two cases when a 5-vertex $v$ is incident to three (5,5,5)-faces.}
\label{fig_lem12}
\end{center}
\end{figure}
\begin{figure}[ht]
\begin{center}
\includegraphics[scale=0.25]{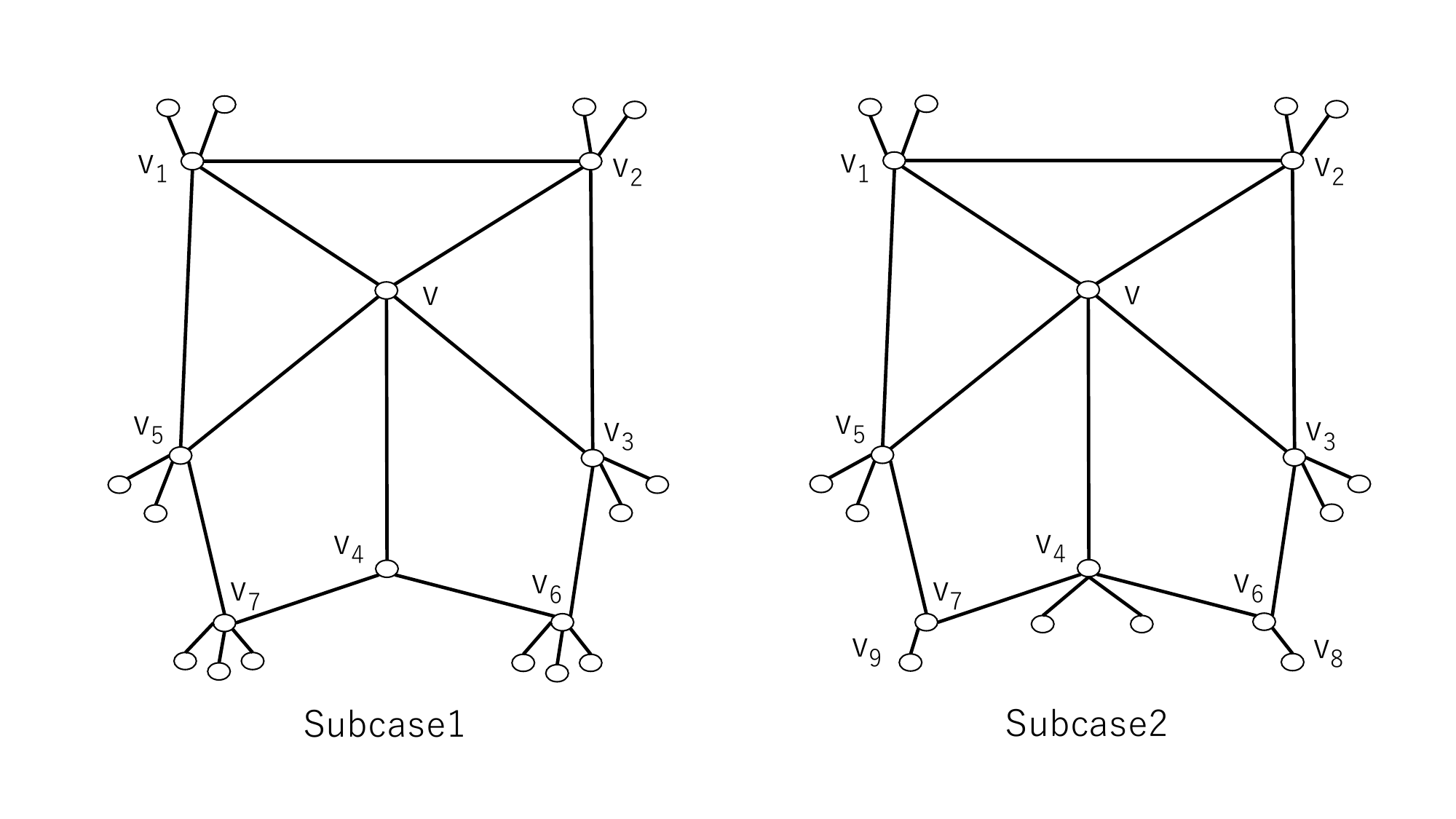}
\caption{Illustrations of Lemma~\ref{lem12} subcases.}
\label{fig_lem12_sub}
\end{center}
\end{figure}

\begin{lem}\label{lem12}
If $v$ is a 5-vertex with $t(v) = 3$ which is incident to three (5,5,5)-faces, then at least one of the other two faces may be (5,5,5,3)-face, but both of them can not be (5,5,5,3)-faces.(See Figure~\ref{fig_lem12}.)
\begin{proof}
Let $N_G(v) = \{v_1,v_2,v_3,v_4,v_5\}$.
Assume $v$ is a 5-vertex which is incident to three (5,5,5)-faces.
We have following two cases.
Case1 : (5,5,5)-face = $[vv_1v_2],[vv_2v_3],[vv_4v_5]$. 
Case2 : (5,5,5)-face = $[vv_1v_2],[vv_2v_3],[vv_5v_1]$.
\begin{itemize}
\item Case1: 
Suppose two (5,5,5,3)-faces are $[vv_5v_6v_1]$ and $[vv_3v_8v_4]$ with $\mathrm{deg}_G(v_6) =\mathrm{deg}_G(v_8) =3$ and $N_G(v_6) = \{v_1,v_5,v_7\}$, $N_G(v_8) = \{v_3,v_4,v_9\}$.
Let $G' = G - \{v,v_6\} + v_1v_5 + v_3v_4 + v_1v_7$.
By the minimality of $G$, $G'$ has a 2-distance 17-coloring $\phi'$.
Let $\phi$ be a coloring of $G$ such that every vertex in $V(G)$, except for $v$ and $v_6$, is colored using $\phi'$.
Since $\Delta(G) \leq 5$, it follows that $|C_\phi(v)| \leq 3 + 3 + 3 + 1 + 3 + 3 = 16$ and $|C| - |C_\phi(v)| \geq 1$, $|C_\phi(v_6)| \leq 5 + 4 + 4 = 13$ and $|C| - |C_\phi(v_6)| \geq 4$.
If $v$ and $v_6$ are colored with $\phi(v) \in C \setminus C_\phi(v)$, $\phi(v_6) \in C \setminus C_\phi(v_6)$ and $\phi(v) \neq \phi(v_6)$, then there exists a coloring $\phi$ of $G$ such that $\chi_2(G) \leq 17$, which is a contradiction.
\item Case2: 
In this case we have following two subcases.(See Figure~\ref{fig_lem12_sub}.)
Subcase1 : The two (5,5,5,3)-faces are $[vv_3v_6v_4]$ and $[vv_4v_7v_5]$ with $\mathrm{deg}_G(v_4) = 3$.
Subcase2 : The two (5,5,5,3)-faces are $[vv_3v_6v_4]$ and $[vv_4v_7v_5]$ with $\mathrm{deg}_G(v_6) = \mathrm{deg}_G(v_7) = 3$ and $N_G(v_6) = \{v_3,v_4,v_8\}$, $N_G(v_7) = \{v_4,v_5,v_9\}$.
\begin{itemize}
\item Subcase1:
Let $G' = G - \{v_4\} + v_6v_7$.
By the minimality of $G$, $G'$ has a 2-distance 17-coloring $\phi'$.
Let $\phi$ be a coloring of $G$ such that every vertex in $V(G)$, except for $v_4$, is colored using $\phi'$.
Since $\Delta(G) \leq 5$, it follows that $|C_\phi(v_4)| \leq 4 + 4 + 5 = 13$ and $|C| - |C_\phi(v_4)| \geq 4$.
If $v_4$ is colored with $\phi(v_4) \in C \setminus C_\phi(v_4)$, then there exists a coloring $\phi$ of $G$ such that $\chi_2(G) \leq 17$, which is a contradiction.
\item Subcase2:
Let $G' = G - \{v,v_7\} + v_3v_4 + v_4v_5 + v_5v_9$.
By the minimality of $G$, $G'$ has a 2-distance 17-coloring $\phi'$.
Let $\phi$ be a coloring of $G$ such that every vertex in $V(G)$, except for $v$ and $v_7$, is colored using $\phi'$.
Since $\Delta(G) \leq 5$, it follows that $|C_\phi(v)| \leq 3 + 3 + 3 + 1 + 3 + 3 = 16$ and $|C| - |C_\phi(v)| \geq 1$, $|C_\phi(v_7)| \leq 5 + 4 + 4 = 13$ and $|C| - |C_\phi(v_7)| \geq 4$.
If $v$ and $v_7$ are colored with $\phi(v) \in C \setminus C_\phi(v)$, $\phi(v_7) \in C \setminus C_\phi(v_7)$ and $\phi(v) \neq \phi(v_7)$, then there exists a coloring $\phi$ of $G$ such that $\chi_2(G) \leq 17$, which is a contradiction.
\end{itemize}
\end{itemize}
\end{proof}
\end{lem}

\begin{figure}[ht]
\begin{center}
\includegraphics[scale=0.25]{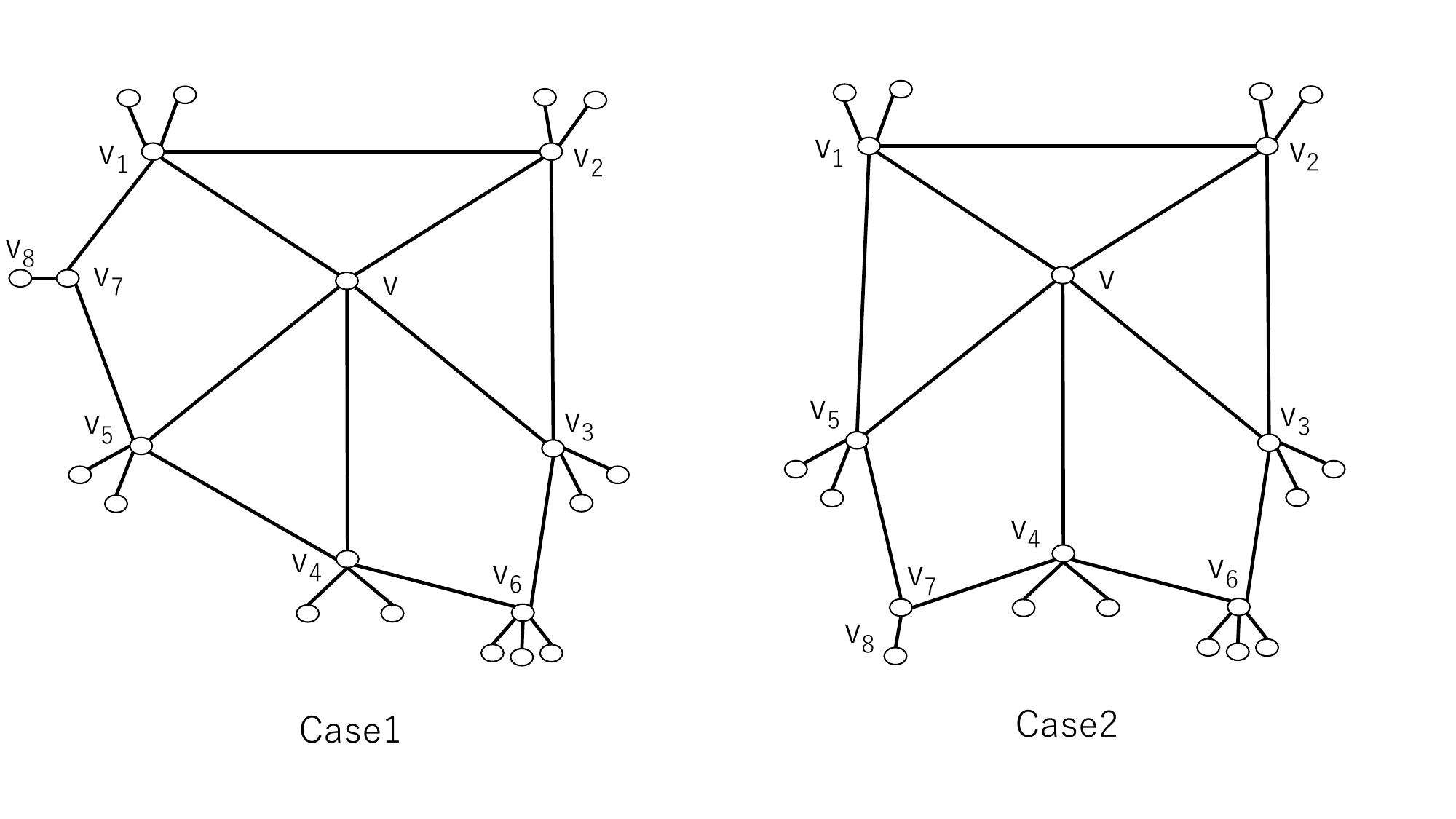}
\caption{Illustrations of Lemma~\ref{lem13}: Assuming a 5-vertex $v$ is incident to three (5,5,5)-faces, one (5,5,5,3)-face and one $(4^+,4^+,4^+,4^+)$-face.}
\label{fig_lem13}
\end{center}
\end{figure}

\begin{lem}\label{lem13}
If $v$ is a 5-vertex with $t(v) = 3$ which is incident to three (5,5,5)-faces and one (5,5,5,3)-face, then the other face can not be $(4^+,4^+,4^+,4^+)$-face.(See Figure~\ref{fig_lem13}.)
\begin{proof}
Let $N_G(v) = \{v_1,v_2,v_3,v_4,v_5\}$.
Assume $v$ is a 5-vertex which is incident to three (5,5,5)-faces and one (5,5,5,3)-face.
We have following two cases.
Case1 : (5,5,5)-face = $[vv_1v_2],[vv_2v_3],[vv_4v_5]$,
(5,5,5,3)-face = $[vv_5v_7v_1]$ with $\mathrm{deg}_G(v_7) = 3$ and $N_G(v_7) = \{v_1,v_5,v_8\}$.
Case2 : (5,5,5)-face = $[vv_1v_2],[vv_2v_3],[vv_5v_1]$,
(5,5,5,3)-face = $[vv_4v_7v_5]$ with $\mathrm{deg}_G(v_7) = 3$ and $N_G(v_7) = \{v_4,v_5,v_8\}$.
\begin{itemize}
\item Case1: 
Suppose $(4^+,4^+,4^+,4^+)$-face is $[vv_3v_6v_4]$.
Let $G' = G - \{v,v_7\} + v_1v_5 + v_3v_4 + v_1v_8$.
By the minimality of $G$, $G'$ has a 2-distance 17-coloring $\phi'$.
Let $\phi$ be a coloring of $G$ such that every vertex in $V(G)$, except for $v$ and $v_7$, is colored using $\phi'$.
Since $\Delta(G) \leq 5$, it follows that $|C_\phi(v)| \leq 3 + 3 + 3 + 1 + 3 + 3 = 16$ and $|C| - |C_\phi(v)| \geq 1$, $|C_\phi(v_7)| \leq 5 + 4 + 4 = 13$ and $|C| - |C_\phi(v_7)| \geq 4$.
If $v$ and $v_7$ are colored with $\phi(v) \in C \setminus C_\phi(v)$, $\phi(v_7) \in C \setminus C_\phi(v_7)$ and $\phi(v) \neq \phi(v_7)$, then there exists a coloring $\phi$ of $G$ such that $\chi_2(G) \leq 17$, which is a contradiction.
\item Case2: 
Suppose $(4^+,4^+,4^+,4^+)$-face is $[vv_3v_6v_4]$.
Let $G' = G - \{v,v_7\} + v_3v_4 + v_4v_5 + v_5v_8$.
By the minimality of $G$, $G'$ has a 2-distance 17-coloring $\phi'$.
Let $\phi$ be a coloring of $G$ such that every vertex in $V(G)$, except for $v$ and $v_7$, is colored using $\phi'$.
Since $\Delta(G) \leq 5$, it follows that $|C_\phi(v)| \leq 3 + 3 + 3 + 1 + 3 + 3 = 16$ and $|C| - |C_\phi(v)| \geq 1$, $|C_\phi(v_7)| \leq 5 + 4 + 4 = 13$ and $|C| - |C_\phi(v_7)| \geq 4$.
If $v$ and $v_7$ are colored with $\phi(v) \in C \setminus C_\phi(v)$, $\phi(v_7) \in C \setminus C_\phi(v_7)$ and $\phi(v) \neq \phi(v_7)$, then there exists a coloring $\phi$ of $G$ such that $\chi_2(G) \leq 17$, which is a contradiction.
\end{itemize}
\end{proof}
\end{lem}

\begin{figure}[ht]
\begin{center}
\includegraphics[scale=0.25]{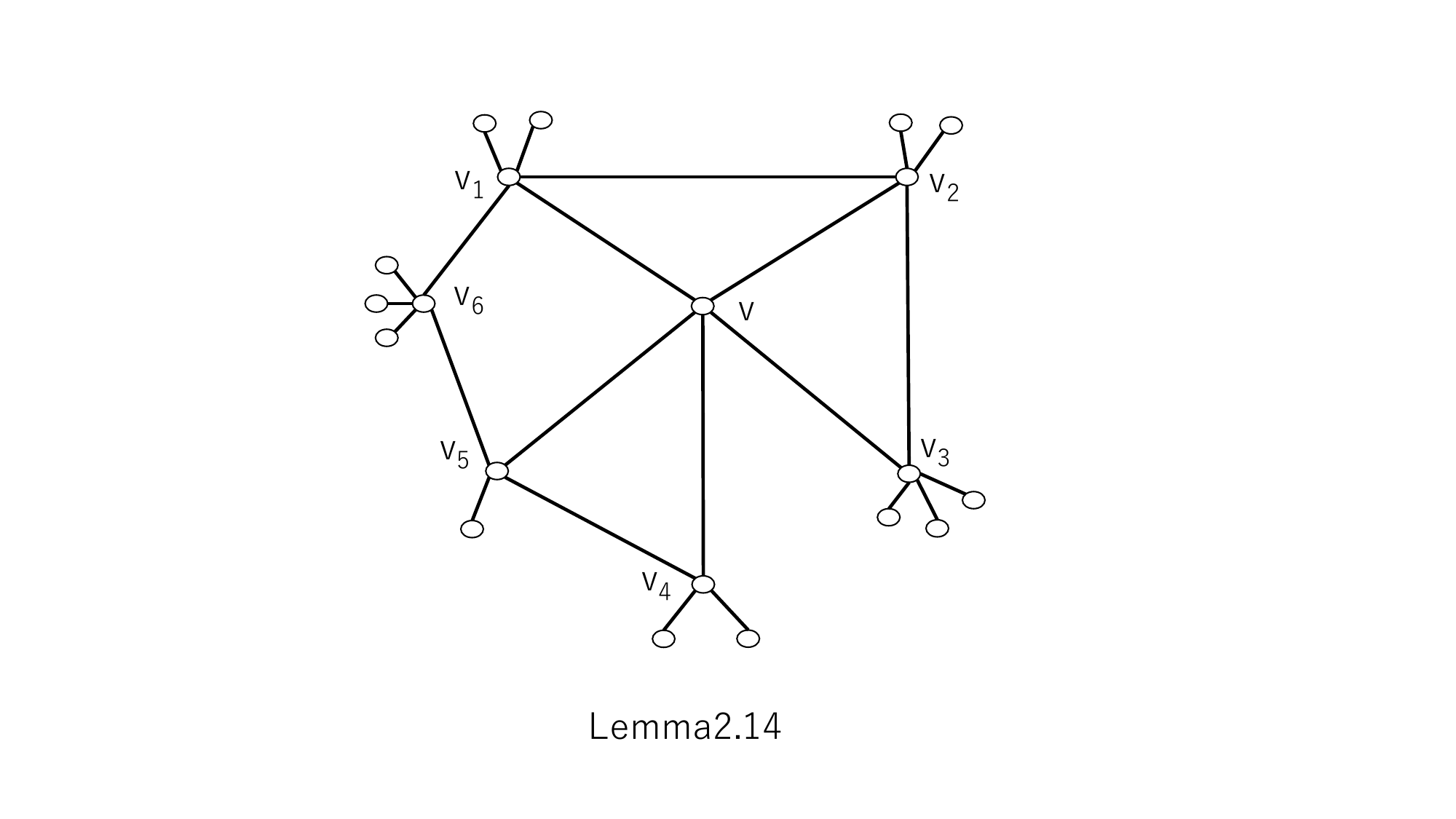}
\caption{Illustration of Lemma~\ref{lem14}: Assuming a 5-vertex $v$ is incident to two (5,5,5)-faces, one (5,4,4)-face and one 4-face.}
\label{fig_lem14}
\end{center}
\end{figure}

\begin{lem}\label{lem14}
If $v$ is a 5-vertex with $t(v) = 3$ which is incident to two (5,5,5)-faces and one (5,4,4)-face, then none of the other faces can be 4-faces.(See Figure~\ref{fig_lem14}.)
\begin{proof}
Let $N_G(v) = \{v_1,v_2,v_3,v_4,v_5\}$.
By Lemma~\ref{lem5} and Lemma~\ref{lem8}, there is no 3-vertex on 3-face and every 4-vertex is incident to at most one 3-face.
Assume $v$ is a 5-vertex which is incident to two (5,5,5)-faces $[vv_1v_2],[vv_2v_3]$ and one (5,4,4)-face $[vv_4v_5]$.
Obviously, $v$ is not incident to (5,5,5,3)-face.
Suppose 4-face is $[vv_5v_6v_1]$.
Let $G' = G - \{v\} + v_1v_5 + v_3v_4$.
By the minimality of $G$, $G'$ has a 2-distance 17-coloring $\phi'$.
Let $\phi$ be a coloring of $G$ such that every vertex in $V(G)$, except for $v$, is colored using $\phi'$.
Since $\Delta(G) \leq 5$, it follows that $|C_\phi(v)| \leq 3 + 3 + 4 + 3 + 2 + 1 = 16$ and $|C| - |C_\phi(v)| \geq 1$.
If $v$ is colored with $\phi(v) \in C \setminus C_\phi(v)$, then there exists a coloring $\phi$ of $G$ such that $\chi_2(G) \leq 17$, which is a contradiction.
\end{proof}
\end{lem}

\begin{figure}[ht]
\begin{center}
\includegraphics[scale=0.25]{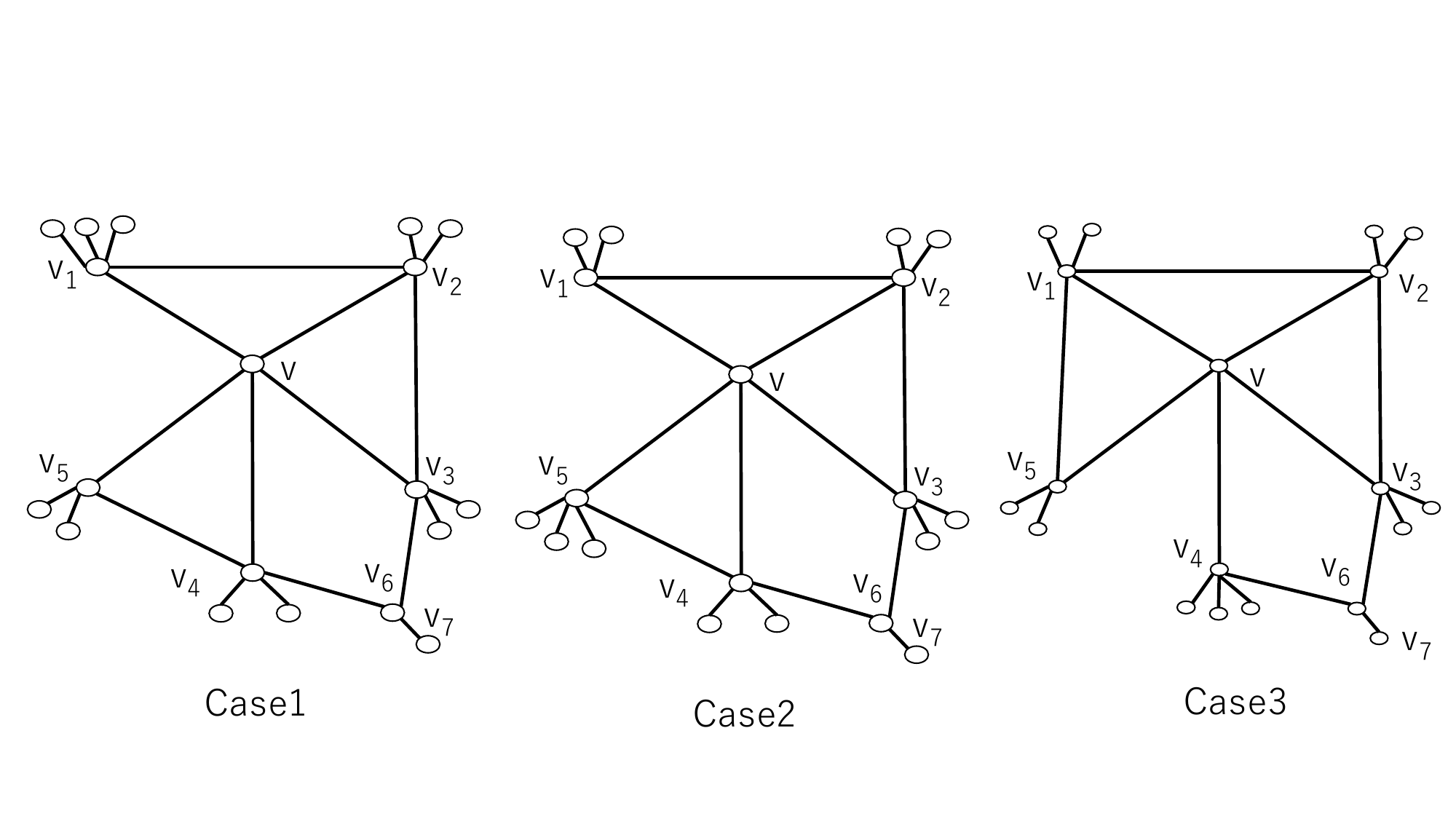}
\caption{Illustrations of Lemma~\ref{lem15}: Assuming a 5-vertex $v$ is incident to two (5,5,5)-faces, one (5,5,4)-face and one (5,5,5,3)-face.}
\label{fig_lem15}
\end{center}
\end{figure}

\begin{lem}\label{lem15}
If $v$ is a 5-vertex with $t(v) = 3$ which is incident to two (5,5,5)-faces and one (5,5,4)-face, then none of the other faces can be (5,5,5,3)-faces.(See Figure~\ref{fig_lem15}.)
\begin{proof}
Let $N_G(v) = \{v_1,v_2,v_3,v_4,v_5\}$.
Assume $v$ is a 5-vertex which is incident to two (5,5,5)-faces and one (5,5,4)-face.
We have following three cases.
Case1 : (5,5,5)-face = $[vv_1v_2],[vv_2v_3]$,
(5,5,4)-face = $[vv_4v_5]$ with $\mathrm{deg}_G(v_5) = 4$.
Case2 : (5,5,5)-face = $[vv_2v_3],[vv_4v_5]$,
(5,5,4)-face = $[vv_1v_2]$ with $\mathrm{deg}_G(v_1) = 4$.
Case3 : (5,5,5)-face = $[vv_1v_2],[vv_2v_3]$,
(5,5,4)-face = $[vv_5v_1]$ with $\mathrm{deg}_G(v_5) = 4$.

\begin{itemize}
\item Case1: 
Suppose (5,5,5,3)-face is $[vv_3v_6v_4]$ with $\mathrm{deg}_G(v_6) = 3$ and $N_G(v_6) = \{v_3,v_4,v_7\}$.
Let $G' = G - \{v,v_6\} + v_1v_5 + v_3v_4 + v_4v_7$.
By the minimality of $G$, $G'$ has a 2-distance 17-coloring $\phi'$.
Let $\phi$ be a coloring of $G$ such that every vertex in $V(G)$, except for $v$ and $v_6$, is colored using $\phi'$.
Since $\Delta(G) \leq 5$, it follows that $|C_\phi(v)| \leq 4 + 3 + 3 + 3 + 3 = 16$ and $|C| - |C_\phi(v)| \geq 1$, $|C_\phi(v_6)| \leq 5 + 4 + 4 = 13$ and $|C| - |C_\phi(v_6)| \geq 4$.
If $v$ and $v_6$ are colored with $\phi(v) \in C \setminus C_\phi(v)$, $\phi(v_6) \in C \setminus C_\phi(v_6)$ and $\phi(v) \neq \phi(v_6)$, then there exists a coloring $\phi$ of $G$ such that $\chi_2(G) \leq 17$, which is a contradiction.
\item Case2: 
Suppose (5,5,5,3)-face is $[vv_3v_6v_4]$ with $\mathrm{deg}_G(v_6) = 3$ and $N_G(v_6) = \{v_3,v_4,v_7\}$.
Let $G' = G - \{v,v_6\} + v_1v_5 + v_3v_4 + v_4v_7$.
By the minimality of $G$, $G'$ has a 2-distance 17-coloring $\phi'$.
Let $\phi$ be a coloring of $G$ such that every vertex in $V(G)$, except for $v$ and $v_6$, is colored using $\phi'$.
Since $\Delta(G) \leq 5$, it follows that $|C_\phi(v)| \leq 3 + 3 + 3 + 3 + 4 = 16$ and $|C| - |C_\phi(v)| \geq 1$, $|C_\phi(v_6)| \leq 5 + 4 + 4 = 13$ and $|C| - |C_\phi(v_6)| \geq 4$.
If $v$ and $v_6$ are colored with $\phi(v) \in C \setminus C_\phi(v)$, $\phi(v_6) \in C \setminus C_\phi(v_6)$ and $\phi(v) \neq \phi(v_6)$, then there exists a coloring $\phi$ of $G$ such that $\chi_2(G) \leq 17$, which is a contradiction.
\item Case3: 
The degree of $v_4$ can not be three.(see Figure~\ref{fig_a1}.)
Suppose (5,5,5,3)-face is $[vv_3v_6v_4]$ with $\mathrm{deg}_G(v_6) = 3$ and $N_G(v_6) = \{v_3,v_4,v_7\}$.
Let $G' = G - \{v,v_6\} + v_4v_5 + v_3v_4 + v_3v_7$.
By the minimality of $G$, $G'$ has a 2-distance 17-coloring $\phi'$.
Let $\phi$ be a coloring of $G$ such that every vertex in $V(G)$, except for $v$ and $v_6$, is colored using $\phi'$.
Since $\Delta(G) \leq 5$, it follows that $|C_\phi(v)| \leq 3 + 3 + 3 + 4 + 3 = 16$ and $|C| - |C_\phi(v)| \geq 1$, $|C_\phi(v_6)| \leq 5 + 4 + 4 = 13$ and $|C| - |C_\phi(v_6)| \geq 4$.
If $v$ and $v_6$ are colored with $\phi(v) \in C \setminus C_\phi(v)$, $\phi(v_6) \in C \setminus C_\phi(v_6)$ and $\phi(v) \neq \phi(v_6)$, then there exists a coloring $\phi$ of $G$ such that $\chi_2(G) \leq 17$, which is a contradiction.
\end{itemize}
\end{proof}
\end{lem}

\begin{figure}[ht]
\begin{center}
\includegraphics[scale=0.25]{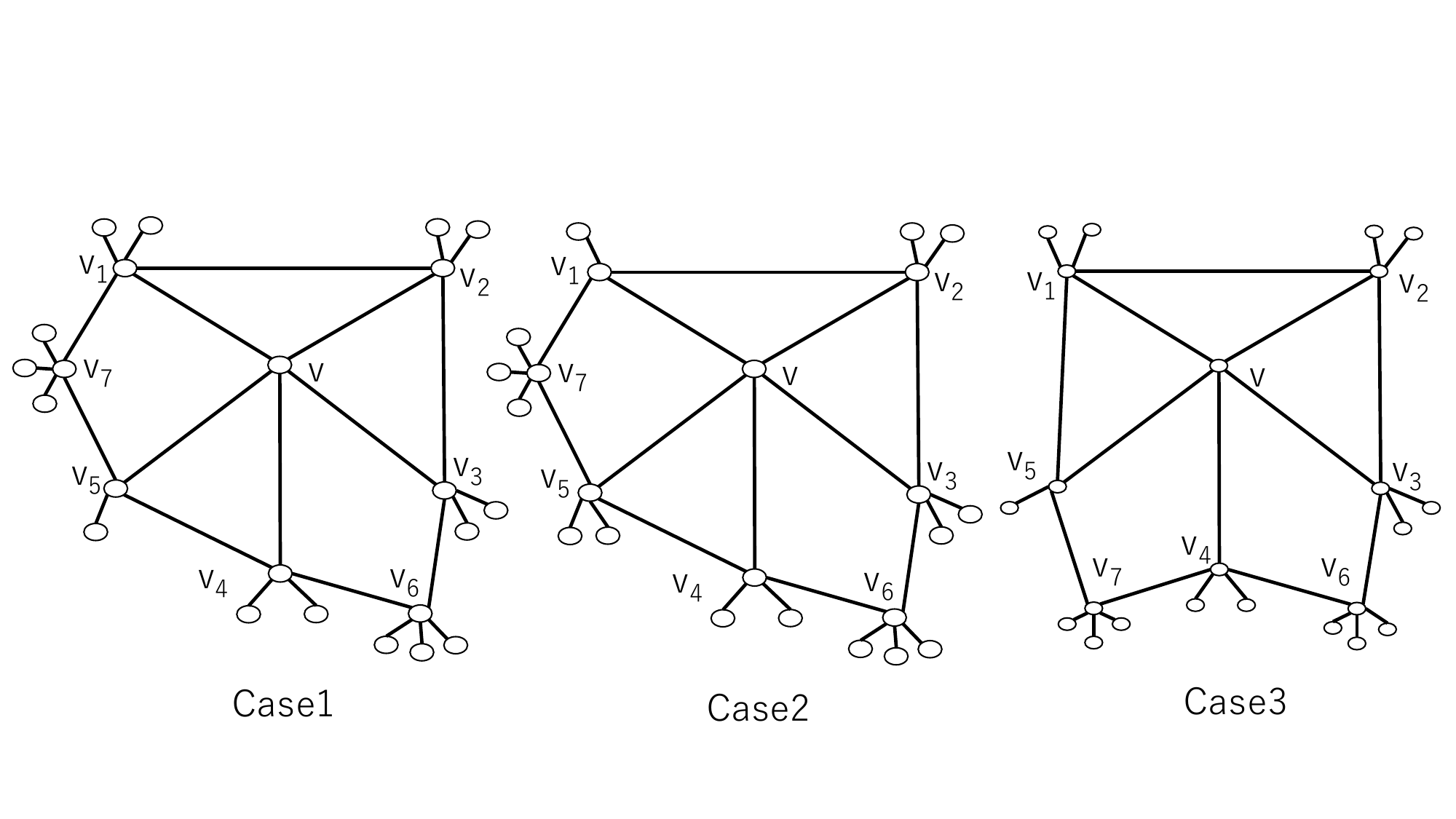}
\caption{Illustrations of Lemma~\ref{lem16}: Assuming a 5-vertex $v$ is incident to two (5,5,5)-faces, one (5,5,4)-face and two $(4^+,4^+,4^+,4^+)$-faces.}
\label{fig_lem16}
\end{center}
\end{figure}

\begin{lem}\label{lem16}
If $v$ is a 5-vertex with t(v) = 3 which is incident to two (5,5,5)-faces and one (5,5,4)-face, then at least one of the other two faces may be $(4^+,4^+,4^+,4^+)$-face, but both of them can not be $(4^+,4^+,4^+,4^+)$-faces.(See Figure~\ref{fig_lem16}.)
\begin{proof}
Let $N_G(v) = \{v_1,v_2,v_3,v_4,v_5\}$.
Assume $v$ is a 5-vertex which is incident to two (5,5,5)-faces and one (5,5,4)-face.
We have following three cases.
Case1 : (5,5,5)-face = $[vv_1v_2],[vv_2v_3]$,
(5,5,4)-face = $[vv_4v_5]$ with $\mathrm{deg}_G(v_5) = 4$.
Case2 : (5,5,5)-face = $[vv_2v_3],[vv_4v_5]$,
(5,5,4)-face = $[vv_1v_2]$ with $\mathrm{deg}_G(v_1) = 4$.
Case3 : (5,5,5)-face = $[vv_1v_2],[vv_2v_3]$,
(5,5,4)-face = $[vv_5v_1]$ with $\mathrm{deg}_G(v_5) = 4$.
\begin{itemize}
\item Case1: 
Suppose two $(4^+,4^+,4^+,4^+)$-faces are $[vv_5v_7v_1]$ and $[vv_3v_6v_4]$.
Let $G' = G - \{v\} + v_1v_5 + v_3v_4$.
By the minimality of $G$, $G'$ has a 2-distance 17-coloring $\phi'$.
Let $\phi$ be a coloring of $G$ such that every vertex in $V(G)$, except for $v$, is colored using $\phi'$.
Since $\Delta(G) \leq 5$, it follows that $|C_\phi(v)| \leq 3 + 3 + 3 + 1 + 3 + 2 + 1 = 16$ and $|C| - |C_\phi(v)| \geq 1$.
If $v$ is colored with $\phi(v) \in C \setminus C_\phi(v)$, then there exists a coloring $\phi$ of $G$ such that $\chi_2(G) \leq 17$, which is a contradiction.
\item Case2: 
Suppose two $(4^+,4^+,4^+,4^+)$-faces are $[vv_5v_7v_1]$ and $[vv_3v_6v_4]$.
Let $G' = G - \{v\} + v_1v_5 + v_3v_4$.
By the minimality of $G$, $G'$ has a 2-distance 17-coloring $\phi'$.
Let $\phi$ be a coloring of $G$ such that every vertex in $V(G)$, except for $v$, is colored using $\phi'$.
Since $\Delta(G) \leq 5$, it follows that $|C_\phi(v)| \leq 2 + 3 + 3 + 1 + 3 + 3 + 1 = 16$ and $|C| - |C_\phi(v)| \geq 1$.
If $v$ is colored with $\phi(v) \in C \setminus C_\phi(v)$, then there exists a coloring $\phi$ of $G$ such that $\chi_2(G) \leq 17$, which is a contradiction.
\item Case3: 
Suppose two $(4^+,4^+,4^+,4^+)$-faces are $[vv_3v_6v_4]$ and $[vv_4v_7v_5]$.
The degree of $v_4$ can not be three.(see Figure~\ref{fig_a1}.)
Let $G' = G - vv_5$.
By the minimality of $G$, $G'$ has a 2-distance 17-coloring $\phi'$.
Every vertex in $V(G)$ is colored using $\phi'$.
Erase the color of $v$ and $v_5$, and recount the available colors for the two vertices.
Since $\Delta(G) \leq 5$, it follows that $|C_{\phi'}(v)| \leq 3 + 3 + 3 + 1 + 3 + 1 + 1 = 15$ and $|C| - |C_{\phi'}(v)| \geq 2$, $|C_{\phi'}(v_5)| \leq 5 + 5 + 3 + 2 = 15$ and $|C| - |C_{\phi'}(v_5)| \geq 2$.
If $v$ and $v_5$ are colored with $\phi'(v) \in C \setminus C_{\phi'}(v)$, $\phi'(v_5) \in C \setminus C_{\phi'}(v_5)$ and $\phi'(v) \neq \phi'(v_5)$, then $\phi'$ can be extended to a 2-distance 17-coloring of $G$, which is a contradiction.
\end{itemize}
\end{proof}
\end{lem}

\begin{figure}[ht]
\begin{minipage}{0.5\hsize}
\begin{center}
\includegraphics[scale=0.2]{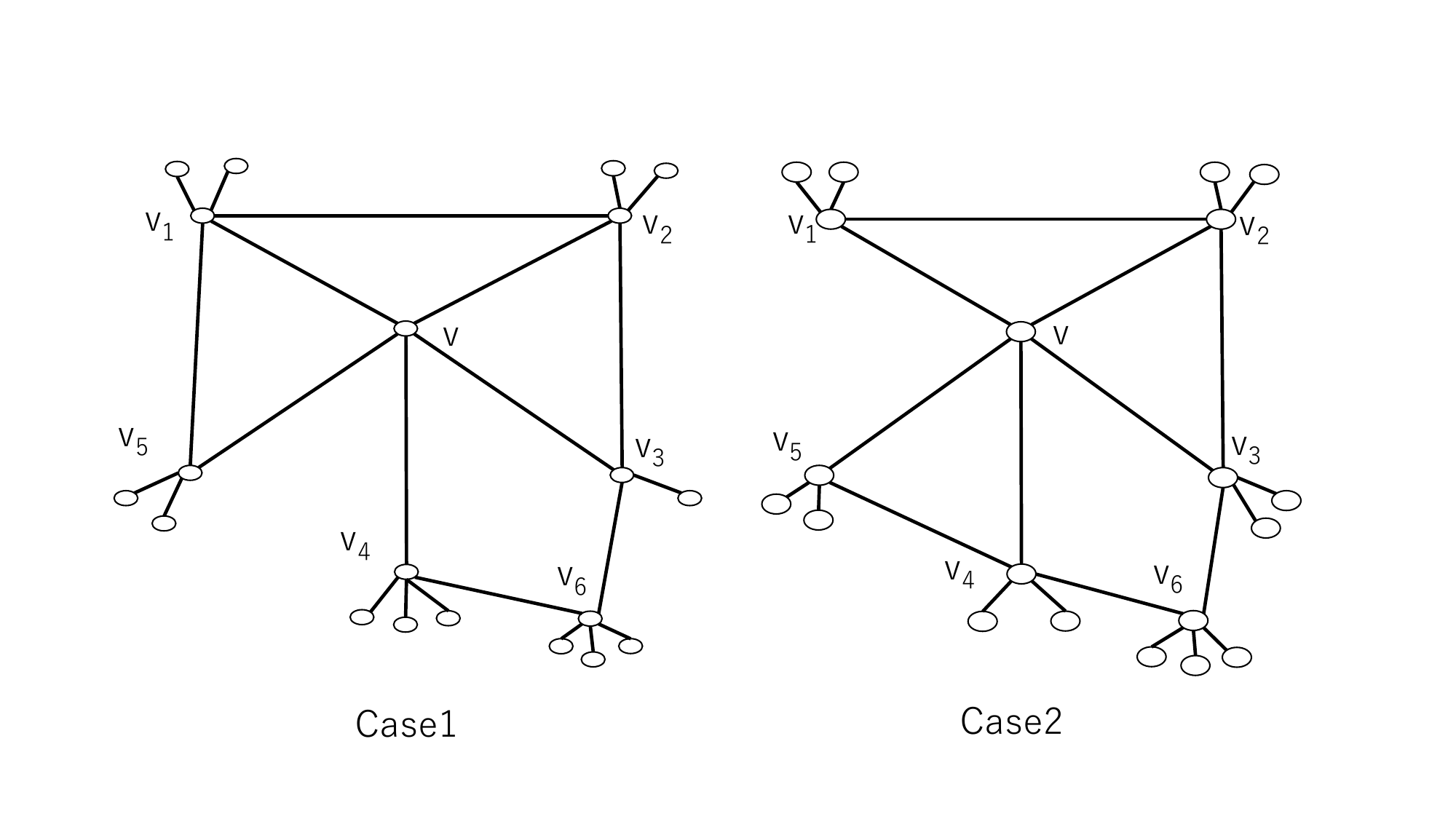}
\label{fig_lem17_1}
\end{center}
\end{minipage}
\begin{minipage}{0.5\hsize}
\begin{center}
\includegraphics[scale=0.2]{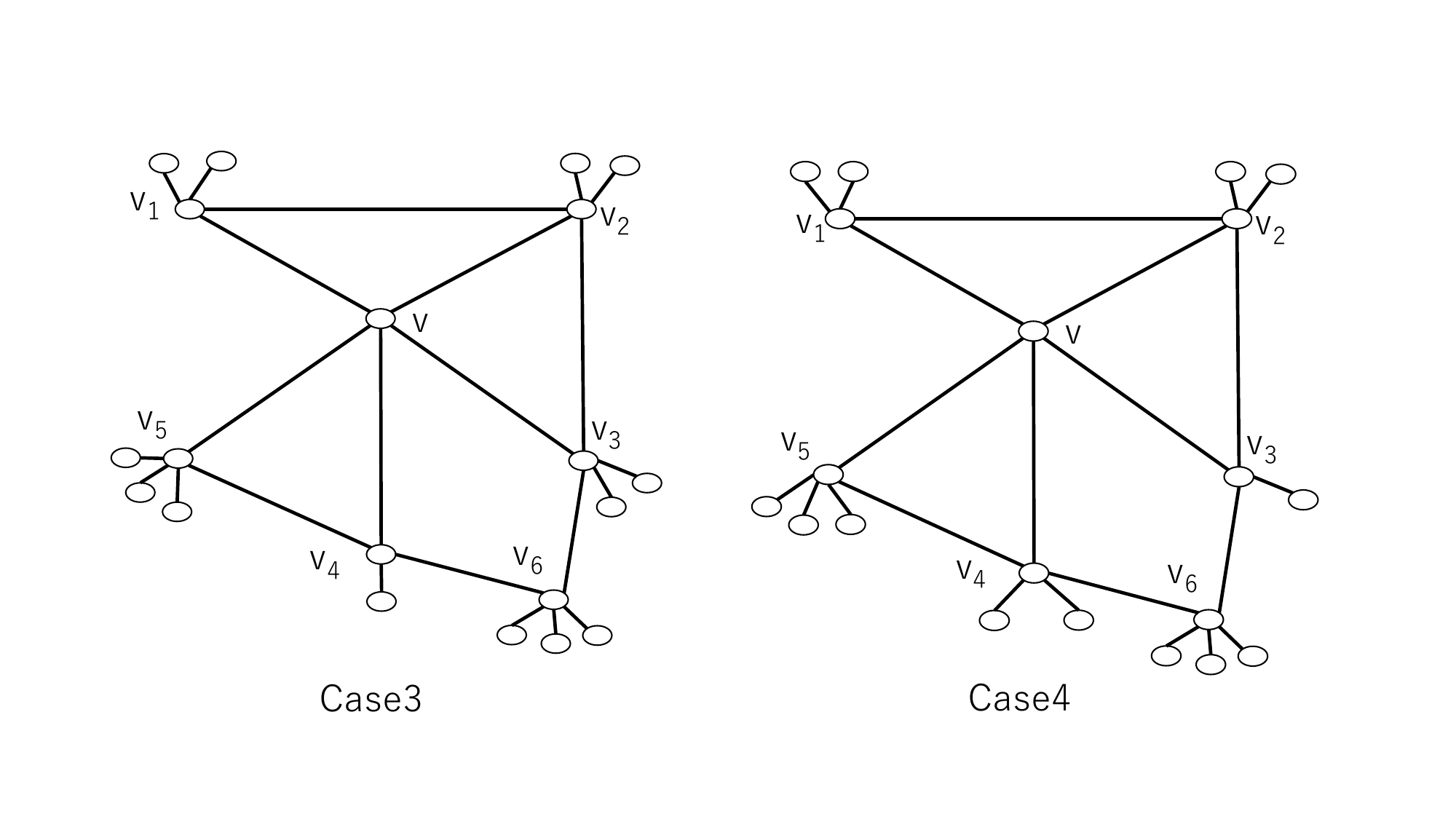}
\label{fig_lem17_2}
\end{center}
\end{minipage}
\caption{Illustrations of Lemma~\ref{lem17}: Assuming a 5-vertex $v$ is incident to one (5,5,5)-face, two (5,5,4)-faces and one 4-face.}
\label{fig_lem17}
\end{figure}

\begin{lem}\label{lem17}
If $v$ is a 5-vertex with $t(v) = 3$ which is incident to one (5,5,5)-face, then none of the other faces can be 4-faces.(See Figure~\ref{fig_lem17}.)
\begin{proof}
Let $N_G(v) = \{v_1,v_2,v_3,v_4,v_5\}$.
Assume $v$ is a 5-vertex which is incident to one (5,5,5)-face.
By Lemma~\ref{lem5}, Lemma~\ref{lem8} and Figure~\ref{fig_a2}, the other two 3-faces incident to $v$ are (5,5,4)-faces.
We have following four cases.
Case1 : (5,5,5)-face = $[vv_1v_2]$,
(5,5,4)-face = $[vv_2v_3],[vv_5v_1]$ with $\mathrm{deg}_G(v_3) = \mathrm{deg}_G(v_5) = 4$.
Case2 : (5,5,5)-face = $[vv_2v_3]$,
(5,5,4)-face = $[vv_1v_2],[vv_4v_5]$ with $\mathrm{deg}_G(v_1) = \mathrm{deg}_G(v_5) = 4$.
Case3 : (5,5,5)-face = $[vv_2v_3]$,
(5,5,4)-face = $[vv_1v_2],[vv_4v_5]$ with $\mathrm{deg}_G(v_1) = \mathrm{deg}_G(v_4) = 4$.
Case4 : (5,5,5)-face = $[vv_4v_5]$,
(5,5,4)-face = $[vv_1v_2],[vv_2v_3]$ with $\mathrm{deg}_G(v_1) = \mathrm{deg}_G(v_3) = 4$.
\begin{itemize}
\item Case1: 
Obviously, $v$ is not incident to any (5,5,5,3)-face.
Suppose $(4^+,4^+,4^+,4^+)$-face is $[vv_3v_6v_4]$.
Let $G' = G - vv_3$.
By the minimality of $G$, $G'$ has a 2-distance 17-coloring $\phi'$.
Every vertex in $V(G)$ is colored using $\phi'$.
Erase the color of $v$ and $v_3$, and recount the available colors for the two vertices.
Since $\Delta(G) \leq 5$, it follows that $|C_{\phi'}(v)| \leq 3 + 3 + 1 + 1 + 4 + 3 = 15$ and $|C| - |C_{\phi'}(v)| \geq 2$, $|C_{\phi'}(v_3)| \leq 5 + 5 + 3 + 2 = 15$ and $|C| - |C_{\phi'}(v_3)| \geq 2$.
If $v$ and $v_3$ are colored with $\phi'(v) \in C \setminus C_{\phi'}(v)$, $\phi'(v_3) \in C \setminus C_{\phi'}(v_3)$ and $\phi'(v) \neq \phi'(v_3)$, then $\phi'$ can be extended to a 2-distance 17-coloring of $G$, which is a contradiction.
\item Case2: 
Suppose $(4^+,4^+,4^+,4^+)$-face is $[vv_3v_6v_4]$.
If $[vv_3v_6v_4]$ is (5,5,5,3)-face, we can prove the same.
Let $G' = G - \{v\} + v_1v_5 + v_3v_4$.
By the minimality of $G$, $G'$ has a 2-distance 17-coloring $\phi'$.
Let $\phi$ be a coloring of $G$ such that every vertex in $V(G)$, except for $v$, is colored using $\phi'$.
Since $\Delta(G) \leq 5$, it follows that $|C_\phi(v)| \leq 3 + 3 + 3 + 1 + 3 + 3 = 16$ and $|C| - |C_\phi(v)| \geq 1$.
If $v$ is colored with $\phi(v) \in C \setminus C_\phi(v)$, then there exists a coloring $\phi$ of $G$ such that $\chi_2(G) \leq 17$, which is a contradiction.
\item Case3: 
Obviously, $v$ is not incident to any (5,5,5,3)-face.
Suppose $(4^+,4^+,4^+,4^+)$-face is $[vv_3v_6v_4]$.
Let $G' = G - \{v\} + v_1v_5 + v_3v_4$.
By the minimality of $G$, $G'$ has a 2-distance 17-coloring $\phi'$.
Let $\phi$ be a coloring of $G$ such that every vertex in $V(G)$, except for $v$, is colored using $\phi'$.
Since $\Delta(G) \leq 5$, it follows that $|C_\phi(v)| \leq 3 + 3 + 3 + 1 + 2 + 4 = 16$ and $|C| - |C_\phi(v)| \geq 1$.
If $v$ is colored with $\phi(v) \in C \setminus C_\phi(v)$, then there exists a coloring $\phi$ of $G$ such that $\chi_2(G) \leq 17$, which is a contradiction.
\item Case4: 
Obviously, $v$ is not incident to any (5,5,5,3)-face.
Suppose $(4^+,4^+,4^+,4^+)$-face is $[vv_3v_6v_4]$.
Let $G' = G - \{v\} + v_1v_5 + v_3v_4$.
By the minimality of $G$, $G'$ has a 2-distance 17-coloring $\phi'$.
Let $\phi$ be a coloring of $G$ such that every vertex in $V(G)$, except for $v$, is colored using $\phi'$.
Since $\Delta(G) \leq 5$, it follows that $|C_\phi(v)| \leq 3 + 3 + 2 + 1 + 3 + 4 = 16$ and $|C| - |C_\phi(v)| \geq 1$.
If $v$ is colored with $\phi(v) \in C \setminus C_\phi(v)$, then there exists a coloring $\phi$ of $G$ such that $\chi_2(G) \leq 17$, which is a contradiction.
\end{itemize}
\end{proof}
\end{lem}

\begin{figure}[ht]
\begin{center}
\includegraphics[scale=0.25]{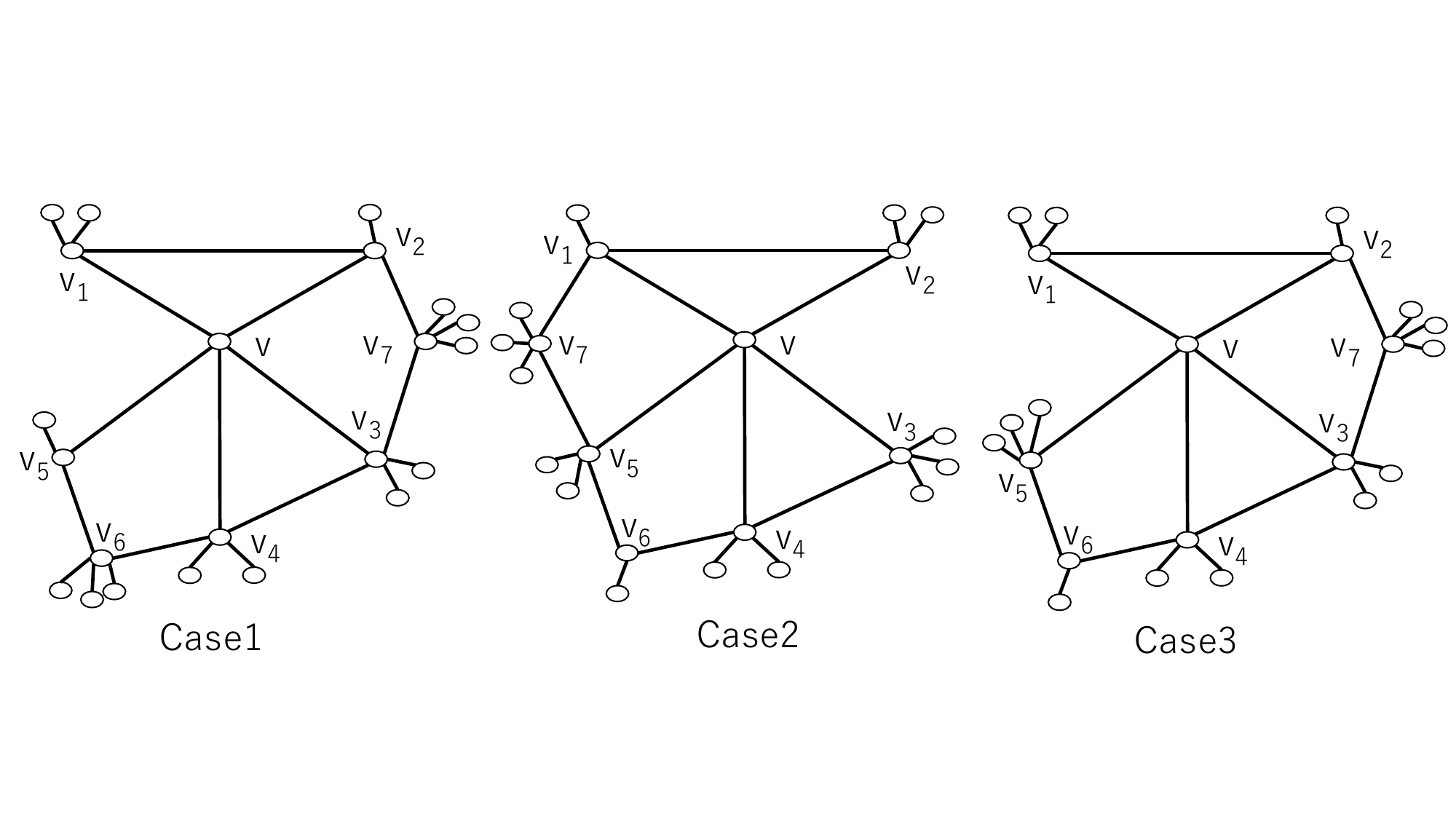}
\end{center}
\begin{center}
\includegraphics[scale=0.25]{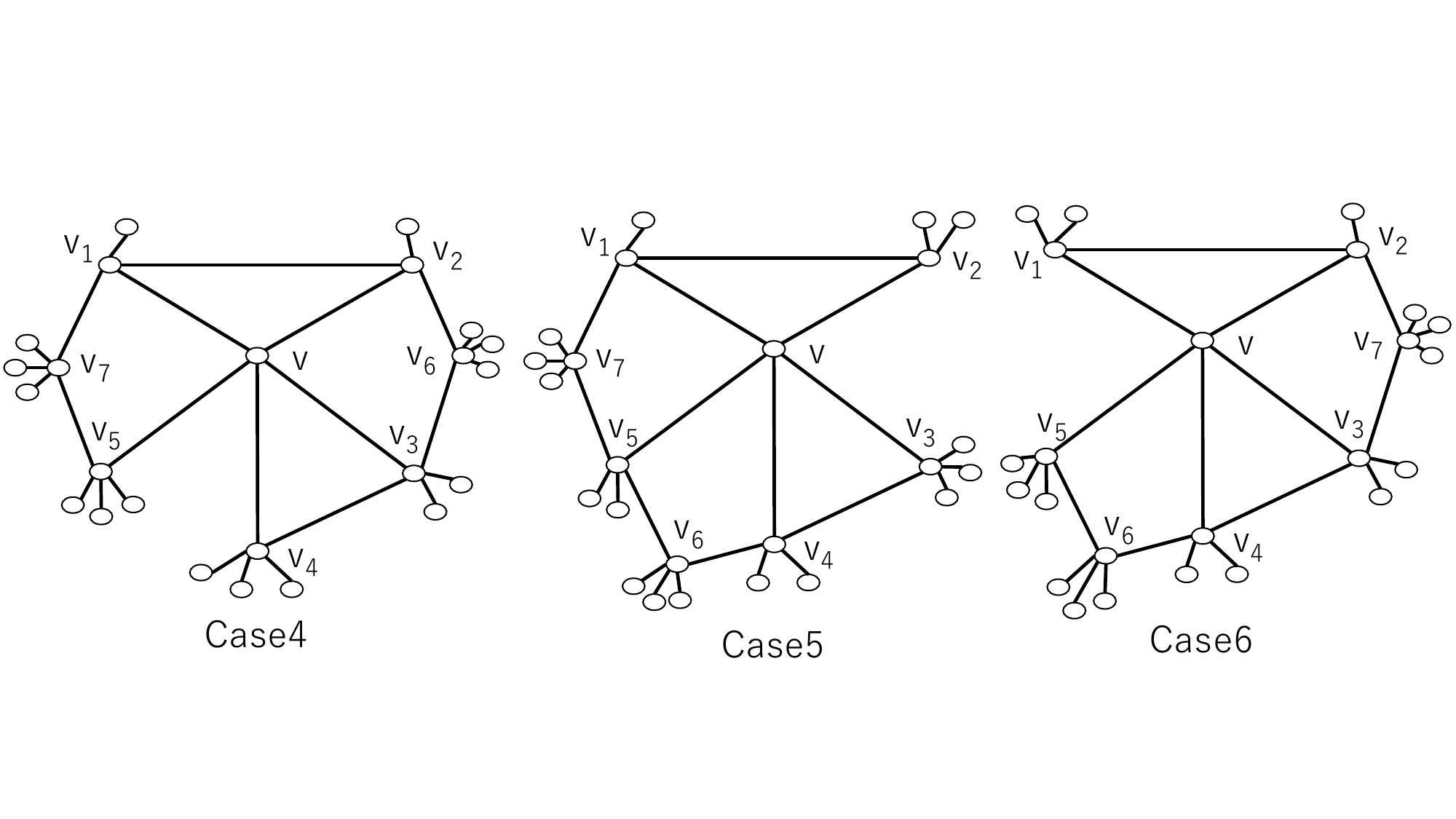}
\end{center}
\caption{Illustrations of Lemma~\ref{lem18}: Assuming a 5-vertex $v$ is incident to one (5,5,5)-face, one (5,5,4)-face and two 4-faces.}
\label{fig_lem18}
\end{figure}

\begin{lem}\label{lem18}
If $v$ is a 5-vertex with $t(v) = 2$ which is incident to one (5,5,5)-face and one (5,4,4)-face, then two of the other three faces can not be 4-faces.(See Figure~\ref{fig_lem18}.)
\begin{proof}
Let $N_G(v) = \{v_1,v_2,v_3,v_4,v_5\}$.
Assume $v$ is a 5-vertex which is incident to one (5,5,5)-face and one (5,4,4)-face.
Let one (5,5,5)-face be $[vv_3v_4]$ and one (5,4,4)-face be $[vv_1v_2]$, we have following six cases.
Case1 : (5,5,5,3)-face = $[vv_4v_6v_5]$ with $\mathrm{deg}_G(v_5) = 3$,
$(4^+,4^+,4^+,4^+)$-face = $[vv_2v_7v_3]$.
Case2 : (5,5,5,3)-face = $[vv_4v_6v_5]$ with $\mathrm{deg}_G(v_6) = 3$,
$(4^+,4^+,4^+,4^+)$-face = $[vv_5v_7v_1]$.
Case3 : (5,5,5,3)-face = $[vv_4v_6v_5]$ with $\mathrm{deg}_G(v_6) = 3$,
$(4^+,4^+,4^+,4^+)$-face = $[vv_2v_7v_3]$.
Case4 : $(4^+,4^+,4^+,4^+)$-face = $[vv_2v_6v_3],[vv_5v_7v_1]$.
Case5 : $(4^+,4^+,4^+,4^+)$-face = $[vv_4v_6v_5],[vv_5v_7v_1]$.
Case6 : $(4^+,4^+,4^+,4^+)$-face = $[vv_2v_7v_3],[vv_4v_6v_5]$.   
\begin{itemize}
\item Case1: Let $G' = G - \{v\} + v_2v_3 + v_4v_5 + v_1v_5$.
By the minimality of $G$, $G'$ has a 2-distance 17-coloring $\phi'$.
Let $\phi$ be a coloring of $G$ such that every vertex in $V(G)$, except for $v$, is colored using $\phi'$.
Since $\Delta(G) \leq 5$, it follows that $|C_\phi(v)| \leq 3 + 2 + 1 + 3 + 3 + 1 + 2 = 15$ and $|C| - |C_\phi(v)| \geq 2$.
If $v$ is colored with $\phi(v) \in C \setminus C_\phi(v)$, then there exists a coloring $\phi$ of $G$ such that $\chi_2(G) \leq 17$, which is a contradiction.
\item Case2: Let $G' = G - vv_1$.
By the minimality of $G$, $G'$ has a 2-distance 17-coloring $\phi'$.
Every vertex in $V(G)$ is colored using $\phi'$.
Erase the color of $v$ and $v_1$, and recount the available colors for the two vertices.
Since $\Delta(G) \leq 5$, it follows that $|C_{\phi'}(v)| \leq 3 + 4 + 3 + 1 + 3 + 1 + 1 = 16$ and $|C| - |C_{\phi'}(v)| \geq 1$, $|C_{\phi'}(v_1)| \leq 5 + 3 + 5 + 1 + 1 = 15$ and $|C| - |C_{\phi'}(v_1)| \geq 2$.
If $v$ and $v_1$ are colored with $\phi'(v) \in C \setminus C_{\phi'}(v)$, $\phi'(v_1) \in C \setminus C_{\phi'}(v_1)$ and $\phi'(v) \neq \phi'(v_1)$, then $\phi'$ can be extended to a 2-distance 17-coloring of $G$, which is a contradiction.
\item Case3: Let $G' = G - vv_2$.
By the minimality of $G$, $G'$ has a 2-distance 17-coloring $\phi'$.
Every vertex in $V(G)$ is colored using $\phi'$.
Erase the color of $v$ and $v_2$, and recount the available colors for the two vertices.
Since $\Delta(G) \leq 5$, it follows that $|C_{\phi'}(v)| \leq 3 + 1 + 1 + 3 + 3 + 1 + 4 = 16$ and $|C| - |C_{\phi'}(v)| \geq 1$, $|C_{\phi'}(v_2)| \leq 5 + 3 + 5 + 1 + 1 = 15$ and $|C| - |C_{\phi'}(v_2)| \geq 2$.
If $v$ and $v_2$ are colored with $\phi'(v) \in C \setminus C_{\phi'}(v)$, $\phi'(v_2) \in C \setminus C_{\phi'}(v_2)$ and $\phi'(v) \neq \phi'(v_2)$, then $\phi'$ can be extended to a 2-distance 17-coloring of $G$, which is a contradiction.
\item Case4: Let $G' = G - vv_1$.
By the minimality of $G$, $G'$ has a 2-distance 17-coloring $\phi'$.
Every vertex in $V(G)$ is colored using $\phi'$.
Erase the color of $v$ and $v_1$, and recount the available colors for the two vertices.
Since $\Delta(G) \leq 5$, it follows that $|C_{\phi'}(v)| \leq 2 + 1 + 3 + 4 + 4 + 1 + 1 = 16$ and $|C| - |C_{\phi'}(v)| \geq 1$, $|C_{\phi'}(v_1)| \leq 5 + 3 + 5 + 1 + 1 = 15$ and $|C| - |C_{\phi'}(v_1)| \geq 2$.
If $v$ and $v_1$ are colored with $\phi'(v) \in C \setminus C_{\phi'}(v)$, $\phi'(v_1) \in C \setminus C_{\phi'}(v_1)$ and $\phi'(v) \neq \phi'(v_1)$, then $\phi'$ can be extended to a 2-distance 17-coloring of $G$, which is a contradiction.
\item Case5: Let $G' = G - vv_1$.
By the minimality of $G$, $G'$ has a 2-distance 17-coloring $\phi'$.
Every vertex in $V(G)$ is colored using $\phi'$.
Erase the color of $v$ and $v_1$, and recount the available colors for the two vertices.
Since $\Delta(G) \leq 5$, it follows that $|C_{\phi'}(v)| \leq 3 + 4 + 3 + 1 + 3 + 1 + 1 = 16$ and $|C| - |C_{\phi'}(v)| \geq 1$, $|C_{\phi'}(v_1)| \leq 5 + 3 + 5 + 1 + 1 = 15$ and $|C| - |C_{\phi'}(v_1)| \geq 2$.
If $v$ and $v_1$ are colored with $\phi'(v) \in C \setminus C_{\phi'}(v)$, $\phi'(v_1) \in C \setminus C_{\phi'}(v_1)$ and $\phi'(v) \neq \phi'(v_1)$, then $\phi'$ can be extended to a 2-distance 17-coloring of $G$, which is a contradiction.
\item Case6: Let $G' = G - vv_2$.
By the minimality of $G$, $G'$ has a 2-distance 17-coloring $\phi'$.
Every vertex in $V(G)$ is colored using $\phi'$.
Erase the color of $v$ and $v_2$, and recount the available colors for the two vertices.
Since $\Delta(G) \leq 5$, it follows that $|C_{\phi'}(v)| \leq 3 + 1 + 1 + 3 + 3 + 1 + 4 = 16$ and $|C| - |C_{\phi'}(v)| \geq 1$, $|C_{\phi'}(v_2)| \leq 5 + 3 + 5 + 1 + 1 = 15$ and $|C| - |C_{\phi'}(v_2)| \geq 2$.
If $v$ and $v_2$ are colored with $\phi'(v) \in C \setminus C_{\phi'}(v)$, $\phi'(v_2) \in C \setminus C_{\phi'}(v_2)$ and $\phi'(v) \neq \phi'(v_2)$, then $\phi'$ can be extended to a 2-distance 17-coloring of $G$, which is a contradiction.
\end{itemize}
\end{proof}
\end{lem}

\begin{figure}[ht]
\begin{center}
\includegraphics[scale=0.25]{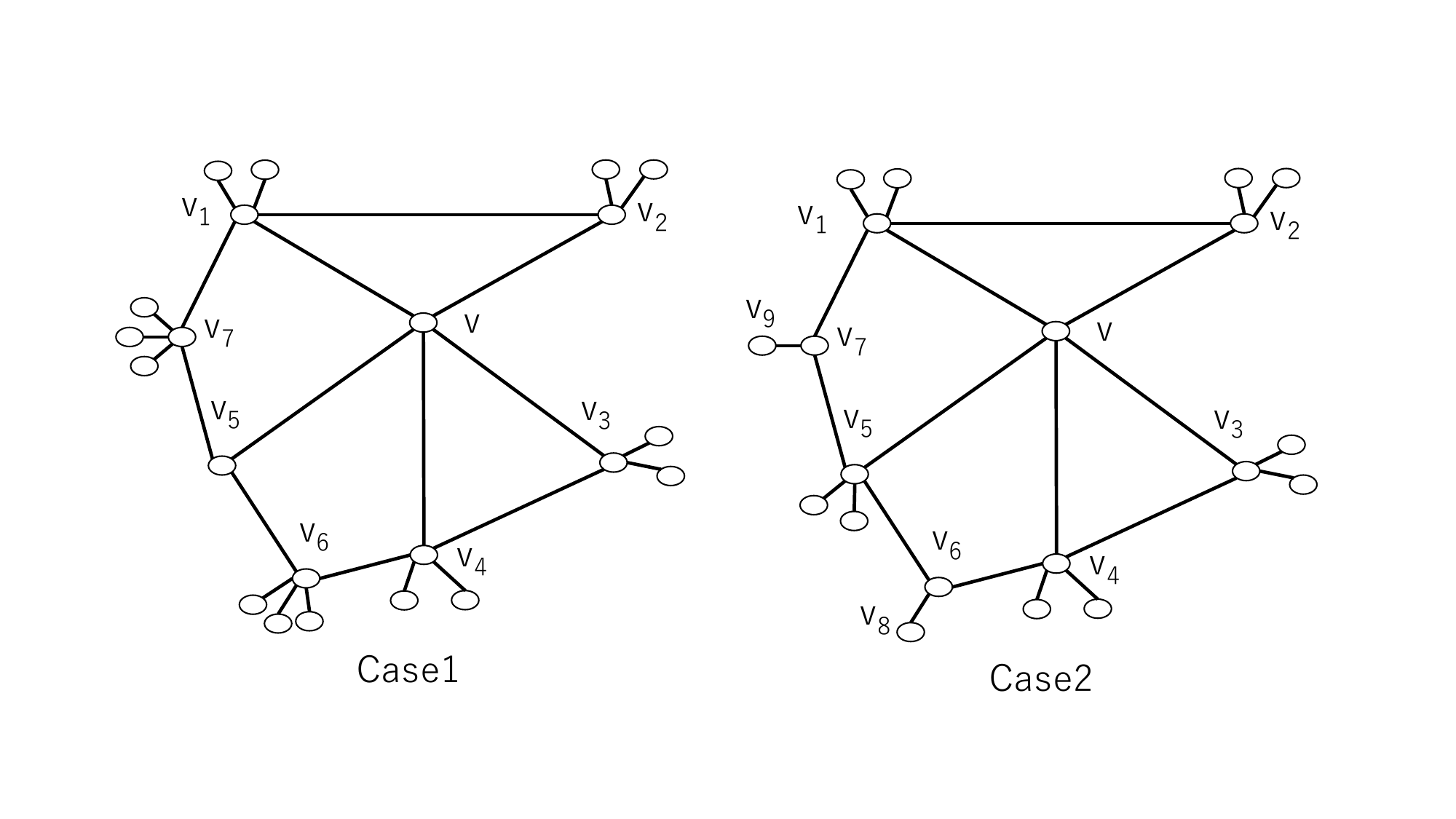}
\caption{Illustrations of Lemma~\ref{lem19}: Assuming a 5-vertex $v$ is incident to two (5,5,4)-faces and two (5,5,5,3)-faces.}
\label{fig_lem19}
\end{center}
\end{figure}

\begin{lem}\label{lem19}
If $v$ is a 5-vertex with $t(v) = 2$ which is incident to two (5,5,4)-faces, then two of the other three faces can not be (5,5,5,3)-faces.(See Figure~\ref{fig_lem19}.)
\begin{proof}
Let $N_G(v) = \{v_1,v_2,v_3,v_4,v_5\}$.
If two (5,5,4)-faces are $[vv_1v_2]$ and $[vv_2v_3]$, then $\mathrm{deg}_G(v_1) = \mathrm{deg}_G(v_3) = 4$ and
$\mathrm{deg}_G(v_2) = 5$ by Lemma~\ref{lem8}. 
In this case, $v$ is not incident to two (5,5,5,3)-face, clearly. 
Assume $v$ is a 5-vertex which is incident to two (5,5,4)-faces $[vv_1v_2]$ and $[vv_3v_4]$. 
Let $\mathrm{deg}_G(v_1) = \mathrm{deg}_G(v_4) = 5$ and $\mathrm{deg}_G(v_2) = \mathrm{deg}_G(v_3) = 4$.
We have following two cases.
Case1 : (5,5,5,3)-face = $[vv_4v_6v_5], [vv_5v_7v_1]$ with $\mathrm{deg}_G(v_5) = 3$.
Case2 : (5,5,5,3)-face = $[vv_4v_6v_5], [vv_5v_7v_1]$ with $\mathrm{deg}_G(v_6) = \mathrm{deg}_G(v_7) = 3$ and $N_G(v_6) = \{v_4,v_5,v_8\}$, $N_G(v_7) = \{v_1,v_5,v_9\}$. 
\begin{itemize}
\item Case1: Let $G' = G - vv_5$.
By the minimality of $G$, $G'$ has a 2-distance 17-coloring $\phi'$.
Every vertex in $V(G)$ is colored using $\phi'$.
Erase the color of $v$ and $v_5$, and recount the available colors for the two vertices.
Since $\Delta(G) \leq 5$, it follows that $|C_{\phi'}(v)| \leq 3 + 3 + 3 + 3 + 1 + 1 = 14$ and $|C| - |C_{\phi'}(v)| \geq 3$, $|C_{\phi'}(v_5)| \leq 5 + 5 + 2 = 12$ and $|C| - |C_{\phi'}(v_5)| \geq 5$.
If $v$ and $v_5$ are colored with $\phi'(v) \in C \setminus C_{\phi'}(v)$, $\phi'(v_5) \in C \setminus C_{\phi'}(v_5)$ and $\phi'(v) \neq \phi'(v_5)$, then $\phi'$ can be extended to a 2-distance 17-coloring of $G$, which is a contradiction.
\item Case2: Let $G' = G - \{v,v_7\} + v_2v_3 + v_4v_5 + v_1v_5 + v_1v_9$.
By the minimality of $G$, $G'$ has a 2-distance 17-coloring $\phi'$.
Let $\phi$ be a coloring of $G$ such that every vertex in $V(G)$, except for $v$ and $v_7$, is colored using $\phi'$.
Since $\Delta(G) \leq 5$, it follows that $|C_\phi(v)| \leq 3 + 3 + 3 + 3 + 1 + 3 = 16 $ and $|C| - |C_\phi(v)| \geq 1$, $|C_\phi(v_7)| \leq 5 + 4 + 4 = 13 $ and $|C| - |C_\phi(v_7)| \geq 4$.
If $v$ and $v_7$ are colored with $\phi(v) \in C \setminus C_\phi(v)$, $\phi(v_7) \in C \setminus C_\phi(v_7)$ and $\phi(v) \neq \phi(v_7)$, then there exists a coloring $\phi$ of $G$ such that $\chi_2(G) \leq 17$, which is a contradiction.
\end{itemize}
\end{proof}
\end{lem}

\begin{figure}[ht]
\begin{center}
\includegraphics[scale=0.25]{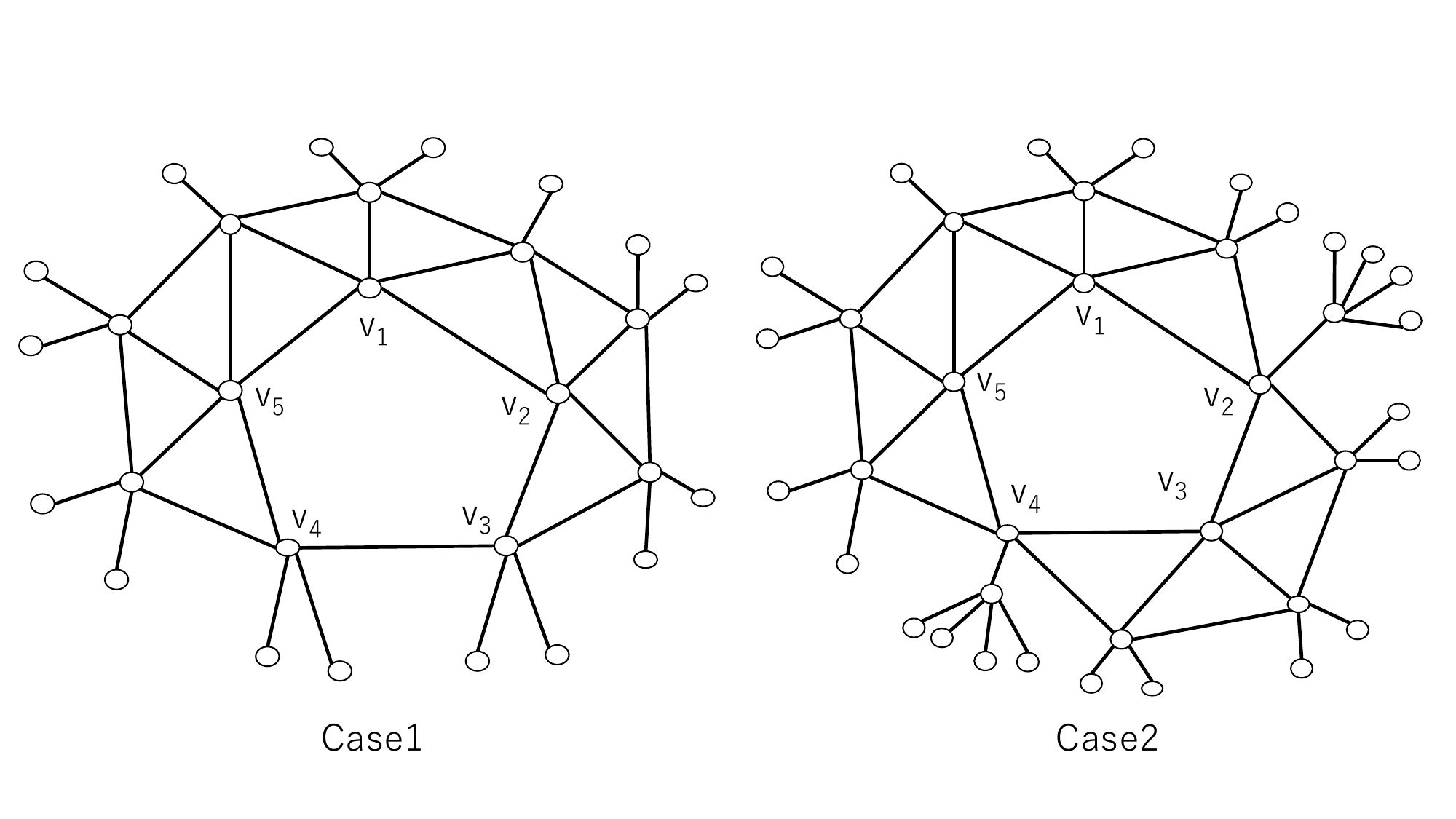}
\caption{Illustrations of Lemma~\ref{lem20}.}
\label{fig_lem20}
\end{center}
\end{figure}

\begin{lem}\label{lem20}
If $[v_1v_2v_3v_4v_5]$ is a 5-face, then there are at most two vertices $v_i$ such that $t(v_i) = 4$ $(i = 1,2,3,4,5)$.(See Figure~\ref{fig_lem20}.)
\begin{proof}
Assume $[v_1v_2v_3v_4v_5]$ is a 5-face and there are three $v_i$ such that $t(v_i) = 4$ $(i = 1,2,3,4,5)$. 
We have following two cases.
Case1 : $t(v_1)=t(v_2)=t(v_5)=4$.
Case2 : $t(v_1)=t(v_3)=t(v_5)=4$.
\begin{itemize}
\item Case1: Let $G' = G - \{v_1\} + v_2v_5$.
By the minimality of $G$, $G'$ has a 2-distance 17-coloring $\phi'$.
Let $\phi$ be a coloring of $G$ such that every vertex in $V(G)$, except for $v_1$, is colored using $\phi'$.
Since $\Delta(G) \leq 5$, it follows that $|C_\phi(v_1)| \leq 3 + 3 + 3 + 3 + 3 = 15$ and $|C| - |C_\phi(v_1)| \geq 2$.
If $v_1$ is colored with $\phi(v_1) \in C \setminus C_\phi(v_1)$, then there exists a coloring $\phi$ of $G$ such that $\chi_2(G) \leq 17$, which is a contradiction.
\item Case2: Let $G' = G - \{v_1\} + v_2v_5$.
By the minimality of $G$, $G'$ has a 2-distance 17-coloring $\phi'$.
Let $\phi$ be a coloring of $G$ such that every vertex in $V(G)$, except for $v_1$, is colored using $\phi'$.
Since $\Delta(G) \leq 5$, it follows that $|C_\phi(v_1)| \leq 3 + 3 + 4 + 3 + 3 = 16$ and $|C| - |C_\phi(v_1)| \geq 1$.
If $v_1$ is colored with $\phi(v_1) \in C \setminus C_\phi(v_1)$, then there exists a coloring $\phi$ of $G$ such that $\chi_2(G) \leq 17$, which is a contradiction.
\end{itemize}
\end{proof}
\end{lem}

\begin{figure}[ht]
\begin{center}
\includegraphics[scale=0.25]{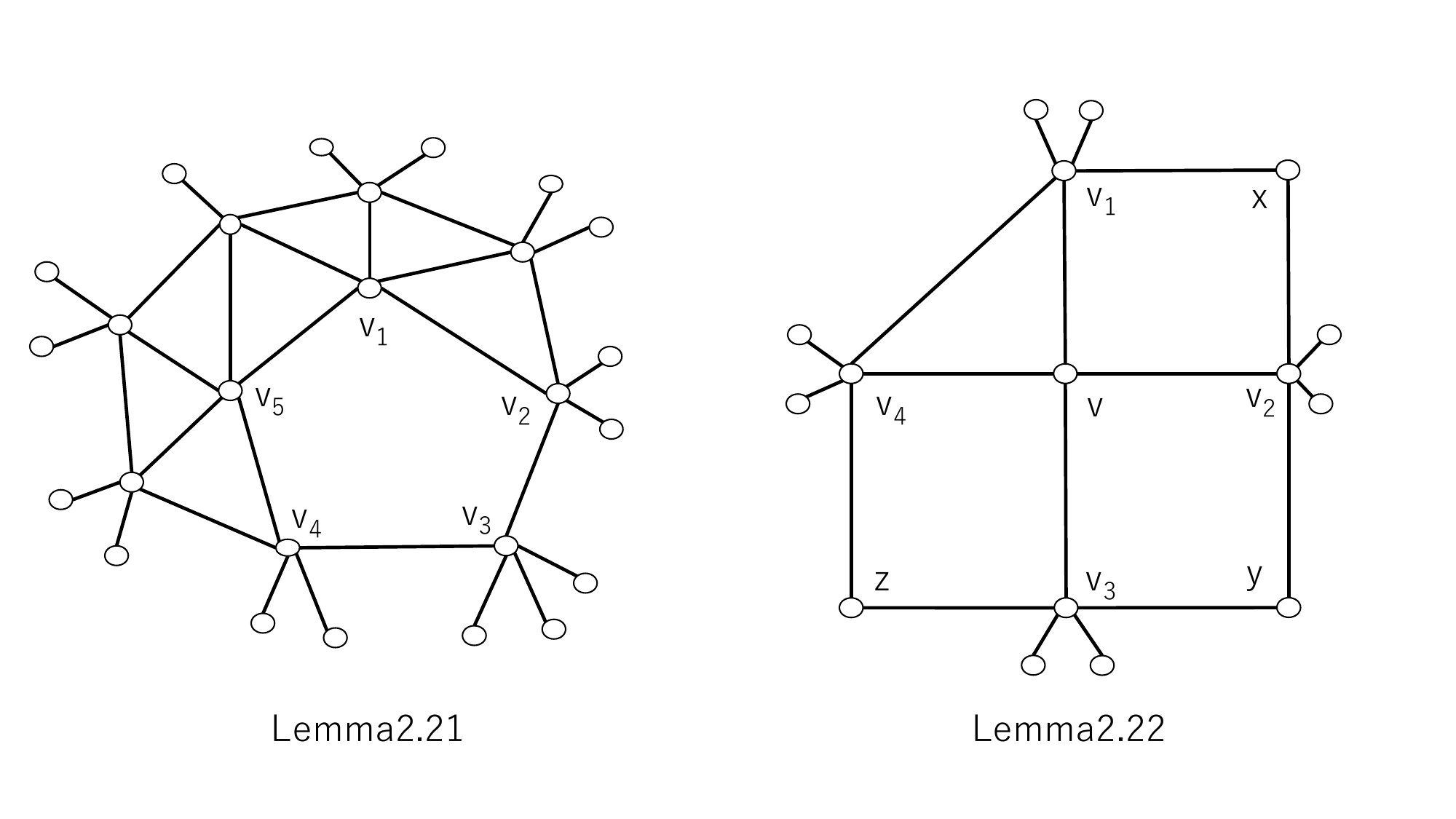}
\caption{Illustrations of Lemma~\ref{lem21} and Lemma~\ref{lem22}.}
\label{fig_lem21_lem22}
\end{center}
\end{figure}

\begin{lem}\label{lem21}
If $[v_1v_2v_3v_4v_5]$ is a 5-face and there exist two vertices $v_i$ such that $t(v_i) = 4$ $(i = 1,2,3,4,5)$, then the two vertices are not adjacent.(See Figure~\ref{fig_lem21_lem22}.)
\begin{proof}
Assume $[v_1v_2v_3v_4v_5]$ is a 5-face with $t(v_1) = t(v_5) = 4$.
Let $G' = G - \{v_1\} + v_2v_5$.
By the minimality of $G$, $G'$ has a 2-distance 17-coloring $\phi'$.
Let $\phi$ be a coloring of $G$ such that every vertex in $V(G)$, except for $v_1$, is colored using $\phi'$.
Since $\Delta(G) \leq 5$, it follows that $|C_\phi(v_1)| \leq 3 + 3 + 4 + 3 + 3 = 16$ and $|C| - |C_\phi(v_1)| \geq 1$.
If $v_1$ is colored with $\phi(v_1) \in C \setminus C_\phi(v_1)$, then there exists a coloring $\phi$ of $G$ such that $\chi_2(G) \leq 17$, which is a contradiction.
\end{proof}
\end{lem}

\begin{lem}\label{lem22}
If $v$ is a 4-vertex with $t(v) = 1$, then the other faces can not be all 4-faces.(See Figure~\ref{fig_lem21_lem22}.)
\begin{proof}
Let $N_G(v) = \{v_1,v_2,v_3,v_4\}$.
Assume $v$ is a 4-vertex which is incident to $[vv_4v_1]$,$[vv_1xv_2]$,$[vv_2yv_3]$ and $[vv_3zv_4]$.
Let $G' = G - \{v\} + v_1v_2 + v_3v_4$.
By the minimality of $G$, $G'$ has a 2-distance 17-coloring $\phi'$.
Let $\phi$ be a coloring of $G$ such that every vertex in $V(G)$, except for $v$, is colored using $\phi'$.
Since $\Delta(G) \leq 5$, it follows that $|C_\phi(v)| \leq 3 + 1 + 3 + 1 + 3 + 1 + 3 = 15$ and $|C| - |C_\phi(v)| \geq 2$.
If $v$ is colored with $\phi(v) \in C \setminus C_\phi(v)$, then there exists a coloring $\phi$ of $G$ such that $\chi_2(G) \leq 17$, which is a contradiction.
\end{proof}
\end{lem}

\begin{figure}[ht]
\begin{minipage}{0.5\hsize}
\begin{center}
\includegraphics[scale=0.2]{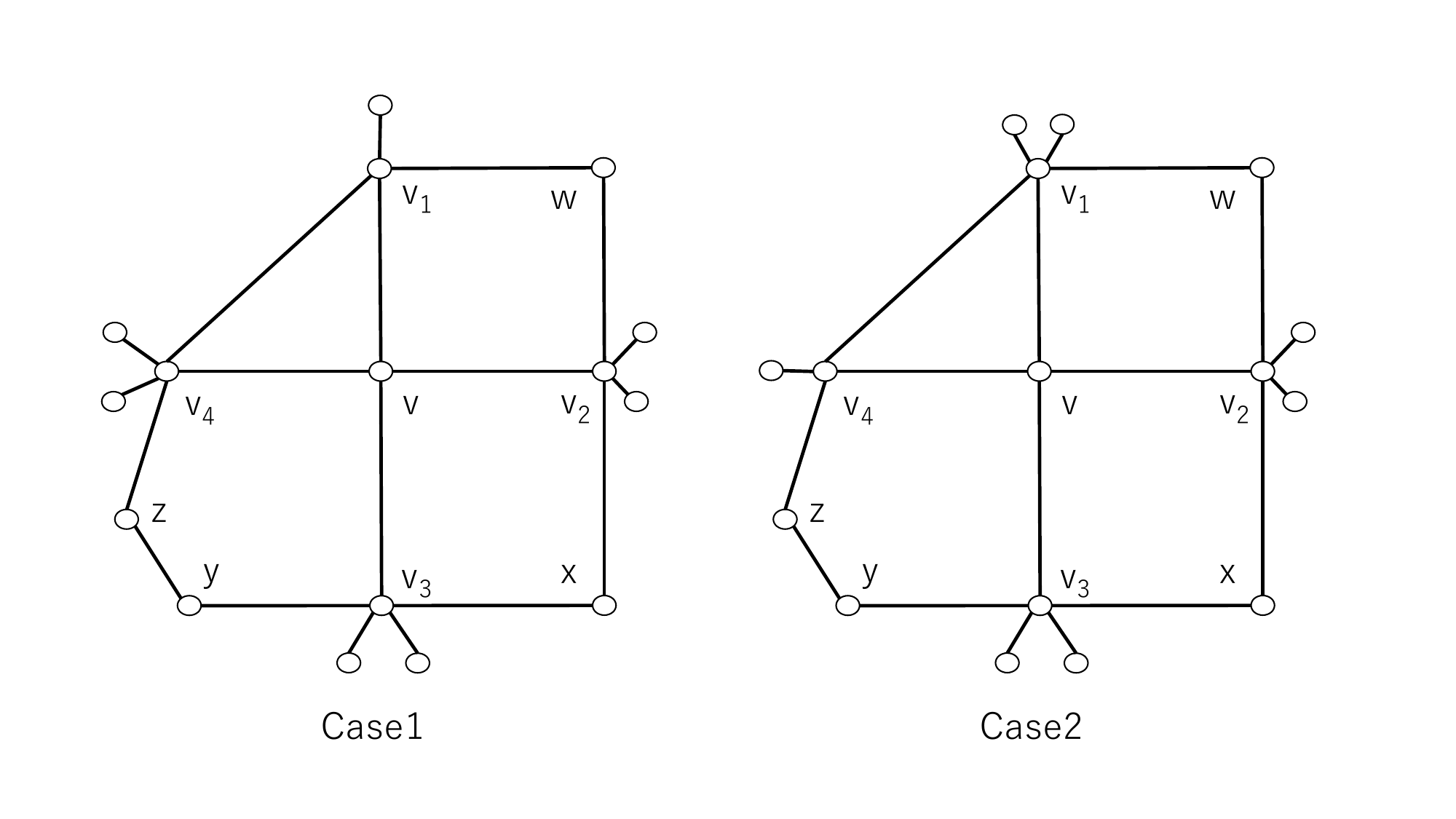}
\label{fig_lem23_1}
\end{center}
\end{minipage}
\begin{minipage}{0.5\hsize}
\begin{center}
\includegraphics[scale=0.2]{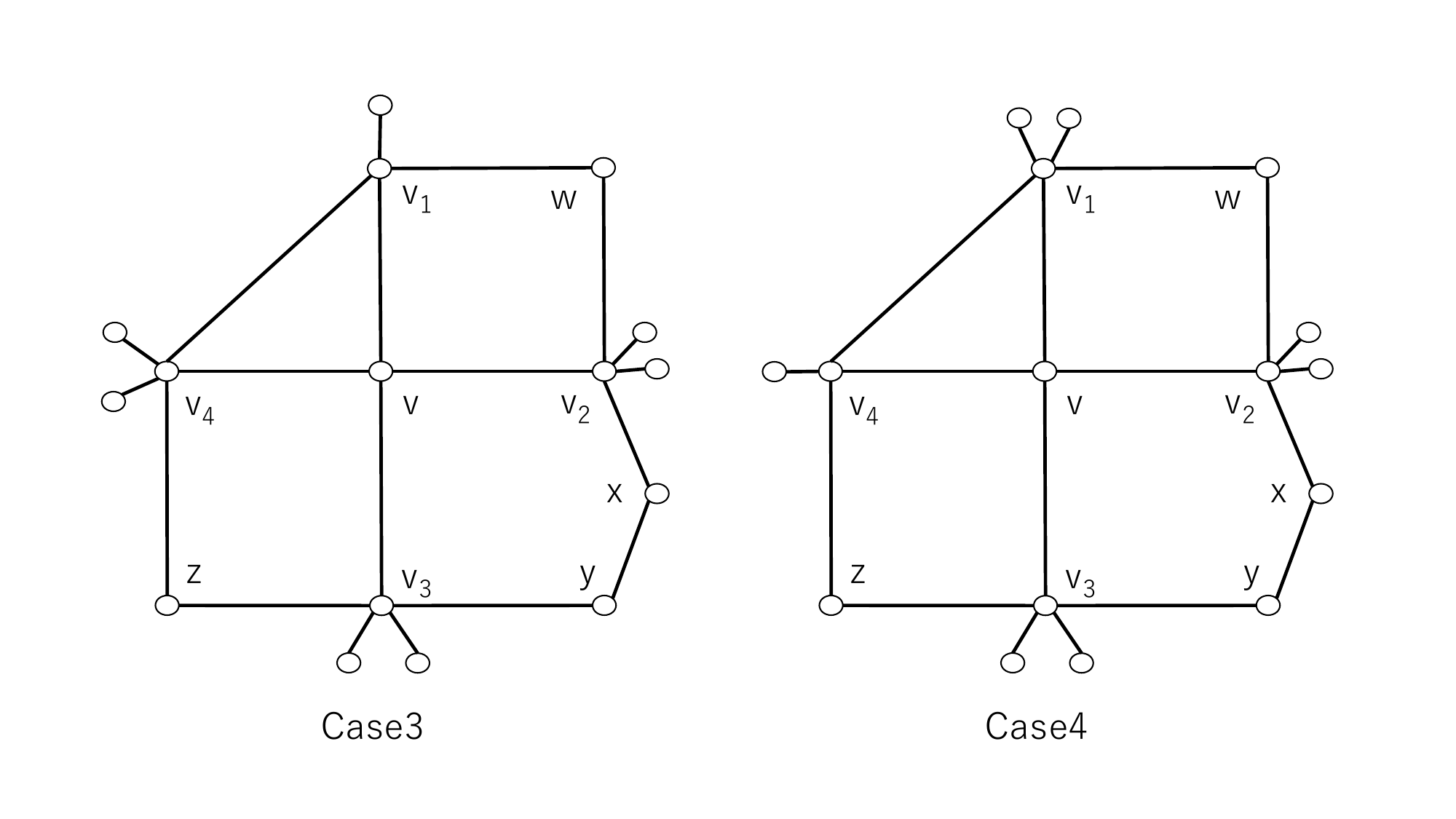}
\label{fig_lem23_2}
\end{center}
\end{minipage}
\caption{Illustrations of Lemma~\ref{lem23}: Assuming a 4-vertex $v$ is incident to one (5,4,4)-face, two 4-faces and one 5-face.}
\label{fig_lem23}
\end{figure}

\begin{lem}\label{lem23}
If $v$ is a 4-vertex with $t(v) = 1$ which is incident to one (5,4,4)-face, then the other faces can not be two 4-faces and one 5-face.(See Figure~\ref{fig_lem23}.)
\begin{proof}
Let $N_G(v) = \{v_1,v_2,v_3,v_4\}$.
Assume $v$ is a 4-vertex which is incident to one (5,4,4)-face $[vv_4v_1]$.
We have following four cases.
Case1 : 4-face = $[vv_1wv_2],[vv_2xv_3]$, 5-face = $[vv_3yzv_4]$ and $\mathrm{deg}_G(v_1) = 4$.
Case2 : 4-face = $[vv_1wv_2],[vv_2xv_3]$, 5-face = $[vv_3yzv_4]$ and $\mathrm{deg}_G(v_4) = 4$.
Case3 : 4-face = $[vv_1wv_2],[vv_3zv_4]$, 5-face = $[vv_2xyv_3]$ and $\mathrm{deg}_G(v_1) = 4$.
Case4 : 4-face = $[vv_1wv_2],[vv_3zv_4]$, 5-face = $[vv_2xyv_3]$ and $\mathrm{deg}_G(v_4) = 4$.
\begin{itemize}
\item Case1: 
Let $G' = G - \{v\} + v_1v_2 + v_3v_4$.
By the minimality of $G$, $G'$ has a 2-distance 17-coloring $\phi'$.
Let $\phi$ be a coloring of $G$ such that every vertex in $V(G)$, except for $v$, is colored using $\phi'$.
Since $\Delta(G) \leq 5$, it follows that $|C_\phi(v)| \leq 2 + 1 + 3 + 1 + 3 + 1 + 1 + 3 = 15$ and $|C| - |C_\phi(v)| \geq 2$.
If $v$ is colored with $\phi(v) \in C \setminus C_\phi(v)$, then there exists a coloring $\phi$ of $G$ such that $\chi_2(G) \leq 17$, which is a contradiction.
\item Case2: 
Let $G' = G - \{v\} + v_1v_2 + v_3v_4$.
By the minimality of $G$, $G'$ has a 2-distance 17-coloring $\phi'$.
Let $\phi$ be a coloring of $G$ such that every vertex in $V(G)$, except for $v$, is colored using $\phi'$.
Since $\Delta(G) \leq 5$, it follows that $|C_\phi(v)| \leq 3 + 1 + 3 + 1 + 3 + 1 + 1 + 2 = 15$ and $|C| - |C_\phi(v)| \geq 2$.
If $v$ is colored with $\phi(v) \in C \setminus C_\phi(v)$, then there exists a coloring $\phi$ of $G$ such that $\chi_2(G) \leq 17$, which is a contradiction.
\item Case3: 
Let $G' = G - \{v\} + v_1v_2 + v_1v_3$.
By the minimality of $G$, $G'$ has a 2-distance 17-coloring $\phi'$.
Let $\phi$ be a coloring of $G$ such that every vertex in $V(G)$, except for $v$, is colored using $\phi'$.
Since $\Delta(G) \leq 5$, it follows that $|C_\phi(v)| \leq 2 + 1 + 3 + 1 + 1 + 3 + 1 + 3 = 15$ and $|C| - |C_\phi(v)| \geq 2$.
If $v$ is colored with $\phi(v) \in C \setminus C_\phi(v)$, then there exists a coloring $\phi$ of $G$ such that $\chi_2(G) \leq 17$, which is a contradiction.
\item Case4: 
Let $G' = G - \{v\} + v_2v_4 + v_3v_4$.
By the minimality of $G$, $G'$ has a 2-distance 17-coloring $\phi'$.
Let $\phi$ be a coloring of $G$ such that every vertex in $V(G)$, except for $v$, is colored using $\phi'$.
Since $\Delta(G) \leq 5$, it follows that $|C_\phi(v)| \leq 3 + 1 + 3 + 1 + 1 + 3 + 1 + 2 = 15$ and $|C| - |C_\phi(v)| \geq 2$.
If $v$ is colored with $\phi(v) \in C \setminus C_\phi(v)$, then there exists a coloring $\phi$ of $G$ such that $\chi_2(G) \leq 17$, which is a contradiction.
\end{itemize}
\end{proof}
\end{lem}

\begin{figure}[ht]
\begin{minipage}{0.5\hsize}
\begin{center}
\includegraphics[scale=0.2]{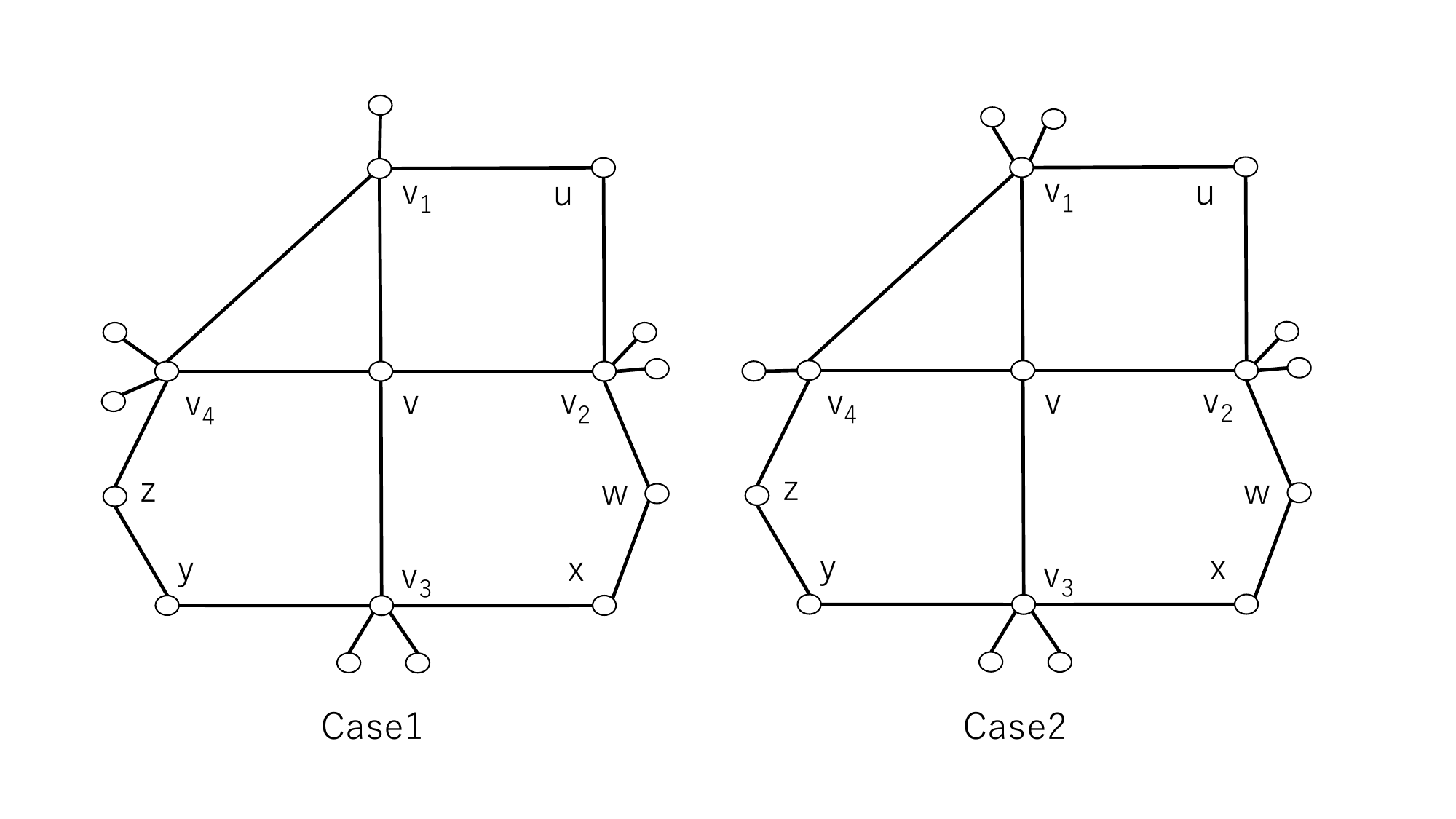}
\label{fig_lem24_1}
\end{center}
\end{minipage}
\begin{minipage}{0.5\hsize}
\begin{center}
\includegraphics[scale=0.2]{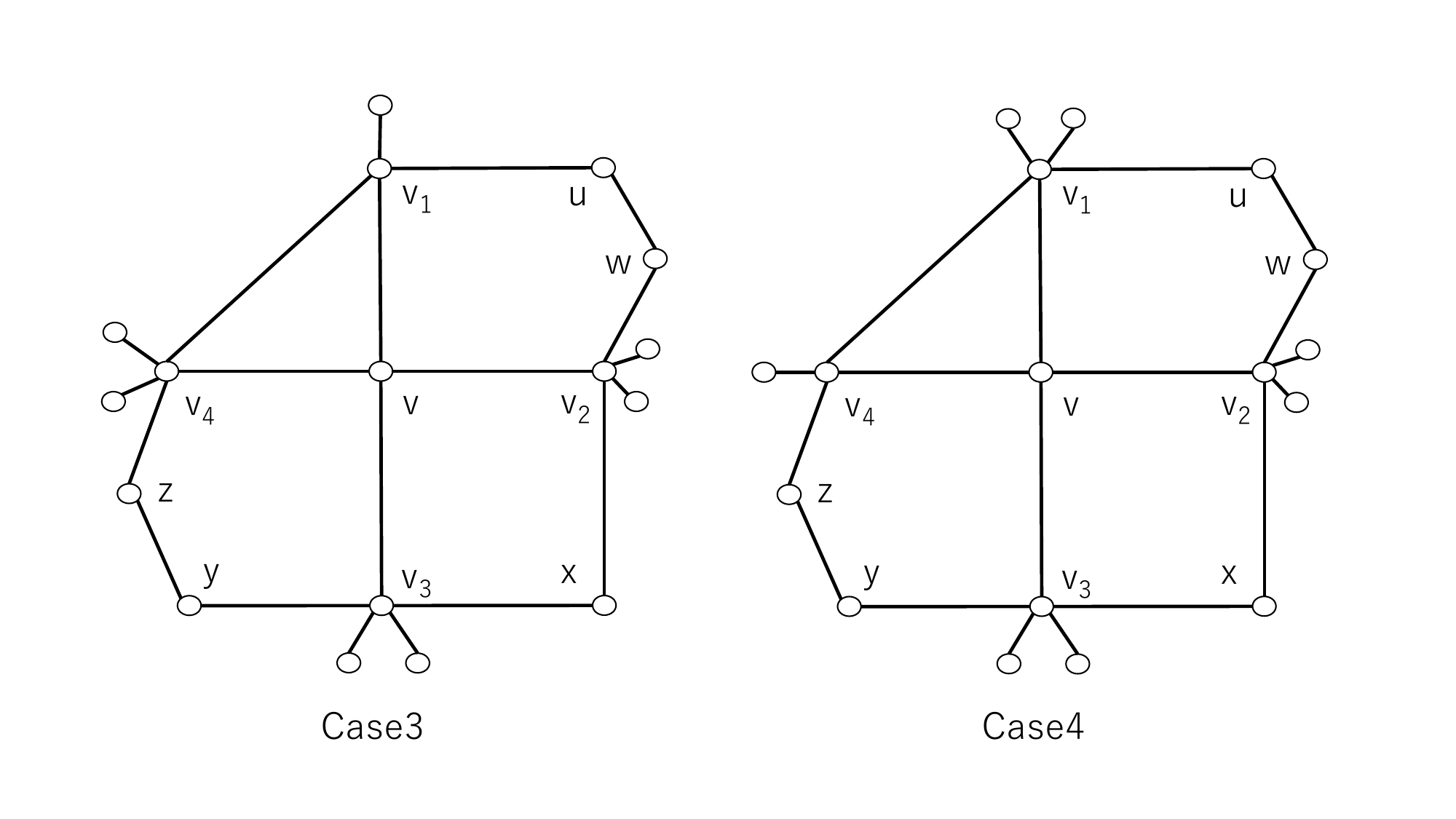}
\label{fig_lem24_2}
\end{center}
\end{minipage}
\caption{Illustrations of Lemma~\ref{lem24}: Assuming a 4-vertex $v$ is incident to one (5,4,4)-face, one 4-face and two 5-faces.}
\label{fig_lem24}
\end{figure}

\begin{lem}\label{lem24}
If $v$ is a 4-vertex with $t(v) = 1$ which is incident to one (5,4,4)-face, then the other faces can not be one 4-face and two 5-faces.(See Figure~\ref{fig_lem24}.)
\begin{proof}
Let $N_G(v) = \{v_1,v_2,v_3,v_4\}$.
Assume $v$ is a 4-vertex which is incident to one (5,4,4)-face $[vv_4v_1]$.
We have following four cases.
Case1 : 4-face = $[vv_1uv_2]$, 5-face = $[vv_2wxv_3],[vv_3yzv_4]$ and $\mathrm{deg}_G(v_1) = 4$.
Case2 : 4-face = $[vv_1uv_2]$, 5-face = $[vv_2wxv_3],[vv_3yzv_4]$ and $\mathrm{deg}_G(v_4) = 4$.
Case3 : 4-face = $[vv_2xv_3]$, 5-face = $[vv_1uwv_2],[vv_3yzv_4]$ and $\mathrm{deg}_G(v_1) = 4$.
Case4 : 4-face = $[vv_2xv_3]$, 5-face = $[vv_1uwv_2],[vv_3yzv_4]$ and $\mathrm{deg}_G(v_4) = 4$.
\begin{itemize}
\item Case1: 
Let $G' = G - \{v\} + v_1v_2 + v_1v_3$.
By the minimality of $G$, $G'$ has a 2-distance 17-coloring $\phi'$.
Let $\phi$ be a coloring of $G$ such that every vertex in $V(G)$, except for $v$, is colored using $\phi'$.
Since $\Delta(G) \leq 5$, it follows that $|C_\phi(v)| \leq 2 + 1 + 3 + 1 + 1 + 3 + 1 + 1 + 3 = 16$ and $|C| - |C_\phi(v)| \geq 1$.
If $v$ is colored with $\phi(v) \in C \setminus C_\phi(v)$, then there exists a coloring $\phi$ of $G$ such that $\chi_2(G) \leq 17$, which is a contradiction.
\item Case2: 
Let $G' = G - \{v\} + v_2v_4 + v_3v_4$.
By the minimality of $G$, $G'$ has a 2-distance 17-coloring $\phi'$.
Let $\phi$ be a coloring of $G$ such that every vertex in $V(G)$, except for $v$, is colored using $\phi'$.
Since $\Delta(G) \leq 5$, it follows that $|C_\phi(v)| \leq 3 + 1 + 3 + 1 + 1 + 3 + 1 + 1 + 2 = 16$ and $|C| - |C_\phi(v)| \geq 1$.
If $v$ is colored with $\phi(v) \in C \setminus C_\phi(v)$, then there exists a coloring $\phi$ of $G$ such that $\chi_2(G) \leq 17$, which is a contradiction.
\item Case3: 
Let $G' = G - \{v\} + v_1v_2 + v_1v_3$.
By the minimality of $G$, $G'$ has a 2-distance 17-coloring $\phi'$.
Let $\phi$ be a coloring of $G$ such that every vertex in $V(G)$, except for $v$, is colored using $\phi'$.
Since $\Delta(G) \leq 5$, it follows that $|C_\phi(v)| \leq 2 + 1 + 1 + 3 + 1 + 3 + 1 + 1 + 3 = 16$ and $|C| - |C_\phi(v)| \geq 1$.
If $v$ is colored with $\phi(v) \in C \setminus C_\phi(v)$, then there exists a coloring $\phi$ of $G$ such that $\chi_2(G) \leq 17$, which is a contradiction.
\item Case4: 
Let $G' = G - \{v\} + v_2v_4 + v_3v_4$.
By the minimality of $G$, $G'$ has a 2-distance 17-coloring $\phi'$.
Let $\phi$ be a coloring of $G$ such that every vertex in $V(G)$, except for $v$, is colored using $\phi'$.
Since $\Delta(G) \leq 5$, it follows that $|C_\phi(v)| \leq 3 + 1 + 1 + 3 + 1 + 3 + 1 + 1 + 2 = 16$ and $|C| - |C_\phi(v)| \geq 1$.
If $v$ is colored with $\phi(v) \in C \setminus C_\phi(v)$, then there exists a coloring $\phi$ of $G$ such that $\chi_2(G) \leq 17$, which is a contradiction.
\end{itemize}
\end{proof}
\end{lem}

\section{Discharging}
In this section, we design discharging rules and complete the proof of Theorem~\ref{main}.
We can derive following equation by Euler's formula $|V(G)| - |E(G)| + |F(G)| = 2$.
$$\sum_{v \in V(G)}(2\mathrm{deg}_G(v)-6) +\sum_{f \in F(G)}(\mathrm{deg}_G(f) - 6) = -12.$$
Note that the total sum is fixed.
We assign an initial charge $2\mathrm{deg}_G(v)-6$ to every vertex and $\mathrm{deg}_G(f)-6$ to every face.
We design appropriate discharging rules and redistribute the charge of the vertices and faces according that rules.
If the final charge of vertices and faces are nonnegative, the following contradiction arises.
$$0 \leq \sum_{v \in V(G)}(2\mathrm{deg}_G(v)-6) +\sum_{f \in F(G)}(\mathrm{deg}_G(f) - 6) = -12 < 0.$$
We design following discharging rules which are improved \cite{chen20222} rules.
\begin{enumerate}
\item[R1] A 5-vertex sends $1$ to each incident (5,5,5)-face.
\item[R2] A 5-vertex sends $\frac{7}{6}$ to each incident (5,5,4)-face.
\item[R3] A 5-vertex sends $1$ to each incident (5,4,4)-face.
\item[R4] A 5-vertex sends $\frac{2}{3}$ to each incident (5,5,5,3)-face.
\item[R5] A 5-vertex sends $\frac{1}{2}$ to each incident $(4^+,4^+,4^+,4^+)$-face.
\item[R6] A 5-vertex sends $\frac{1}{3}$ to each incident 5-face.
\item[R7] A 4-vertex sends $\frac{2}{3}$ to each incident (5,5,4)-face.
\item[R8] A 4-vertex sends $1$ to each incident (5,4,4)-face.
\item[R9] A 4-vertex sends $\frac{1}{2}$ to each incident 4-face.
\item[R10] A 4-vertex sends $\frac{1}{3}$ to each incident 5-face.
\item[R11] A 5-vertex with $t(v) = 4$ does not sends charge except for (5,5,5)-face.
\end{enumerate}

Next, we check the final charge of vertices and faces.
Obviously, the final charge of 3-vertex and $6^+$-face are nonnegative.
Thus, we only check the final charge of 4-vertex, 5-vertex, 3-face, 4-face and 5-face.
First, we show that the final charge for each face is nonnegative.

\begin{enumerate}
\item[3-face:] By Lemma~\ref{lem5} and Lemma~\ref{lem6}, each 3-face is a (5,5,5)-face, a (5,5,4)-face or a (5,4,4)-face.
By R1, the final charge of (5,5,5)-face is $-3 + 3 \times 1 = 0$.
By R2 and R7, the final charge of (5,5,4)-face is $-3 + 2 \times \frac{7}{6} + \frac{2}{3} = 0$.
By R3 and R8, the final charge of (5,4,4)-face is $-3 + 1 +  2 \times 1 = 0$.
\item[4-face:] By Lemma~\ref{lem9}, each 4-face is a (5,5,5,3)-face or a $(4^+,4^+,4^+,4^+)$-face.
By R4, the final charge of (5,5,5,3)-face is $-2 + 3 \times \frac{2}{3} = 0$. 
By R5 and R9, the final charge of $(4^+,4^+,4^+,4^+)$-face is $-2 + 4 \times \frac{1}{2} = 0$. 
\item[5-face:]
By Lemma~\ref{lem4} and Lemma~\ref{lem7}, we have eleven cases about a 5-face.(see Figure~\ref{fig_5-face}.)
The number in the figure indicate the degree of the vertex.
There are at most two vertices $v$ with $t(v) = 4$ on a five face from Lemma~\ref{lem20}.
By Lemma~\ref{lem21} if there are two vertices $v$ with $t(v) = 4$ on a five face, then the two vertices are not adjacent.
A 5-vertex $v$ with $t(v) = 4$ is not adjacent to $4^-$-vertex by Lemma~\ref{lem11}.
\begin{figure}[ht]
\begin{center}
\includegraphics[scale=0.25]{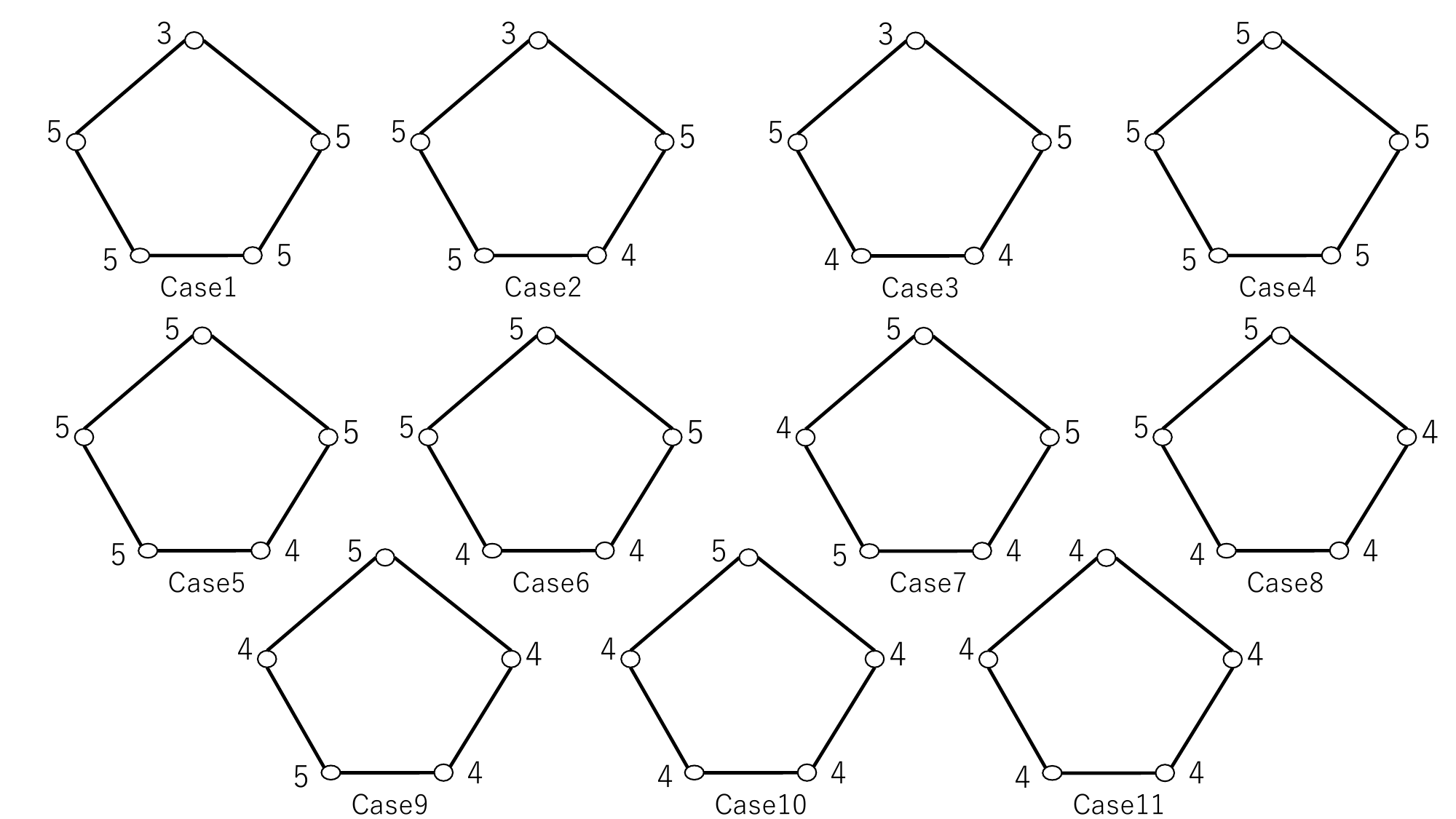}
\caption{Illustrations of eleven cases about  5-faces.}
\label{fig_5-face}
\end{center}
\end{figure}
In Case8 to Case11, the final charge is nonnegative by R10.
Case1: There are at most one 5-vertex with $t(v)=4$.
By R6 and R11, the final charge is at least $-1 + 3 \times \frac{1}{3} =  0$.
Case2: By R6 and R10, the final charge is $-1 + 4 \times \frac{1}{3} = \frac{1}{3}$.
Case3: By R6 and R10, the final charge is $-1 + 4 \times \frac{1}{3} = \frac{1}{3}$.
Case4: There are at most two 5-vertices with $t(v)$ = 4.
By R6 and R11, the final charge is at least $-1 + 3 \times \frac{1}{3} =  0$.
Case5: There are at most one 5-vertex with $t(v)=4$.
By R6, R10 and R11, the final charge is at least $-1 + 4 \times \frac{1}{3} = \frac{1}{3}$.
Case6: There are at most one 5-vertex with $t(v)=4$.
By R6, R10 and R11, the final charge is at least $-1 + 4 \times \frac{1}{3} =  \frac{1}{3}$.
Case7: By R6 and R10, the final charge is $-1 + 5 \times \frac{1}{3} = \frac{2}{3}$.
\end{enumerate}
Next, we show that the final charge for each vertex is nonnegative.
\begin{enumerate}
\item[4-vertex:] By Lemma~\ref{lem8}, every 4-vertex is incident to at most one 3-face.
\begin{itemize}
\item The case where a 4-vertex $v$ is incident to one 3-face.

If a 4-vertex $v$ is incident to one 3-face, then we have following two cases.
\begin{itemize}
\item The case where a 4-vertex $v$ is incident to a (5,5,4)-face.

\quad By Lemma~\ref{lem22}, the 4-vertex $v$ is not incident to three 4-faces.
If the 4-vertex $v$ is incident to a (5,5,4)-face, then we have following three cases.
Case1: The 4-vertex $v$ is incident to two 4-faces and one 5-face.
By R7, R9 and R10, the final charge is $2 - \frac{2}{3} - 2 \times \frac{1}{2} - \frac{1}{3}  = 0$.
Case2: The 4-vertex $v$ is incident to one 4-face and two 5-faces.
By R7, R9 and R10, the final charge is $2 - \frac{2}{3} - \frac{1}{2} - 2 \times \frac{1}{3} = \frac{1}{6}$.
Case3: The 4-vertex $v$ is incident to three 5-faces.
By R7 and R10, the final charge is $2 - \frac{2}{3} - 3 \times \frac{1}{3} = \frac{1}{3}$.
\item The case where a 4-vertex $v$ is incident to (5,4,4)-face.

\quad By Lemma~\ref{lem22}, Lemma~\ref{lem23} and Lemma~\ref{lem24}, the 4-vertex $v$ is not incident to any 4-face.
The 4-vertex $v$ is incident to three 5-faces. 
By R8 and R10, the final charge is $2 - 1 - 3 \times \frac{1}{3} = 0$.
\end{itemize}
\item The case where a 4-vertex $v$ is not incident to any 3-face.

If a 4-vertex $v$ is not incident to any 3-face, then we have five cases.(see Figure~\ref{fig_4-face}.)
\begin{figure}[t]
\begin{center}
\includegraphics[scale=0.25]{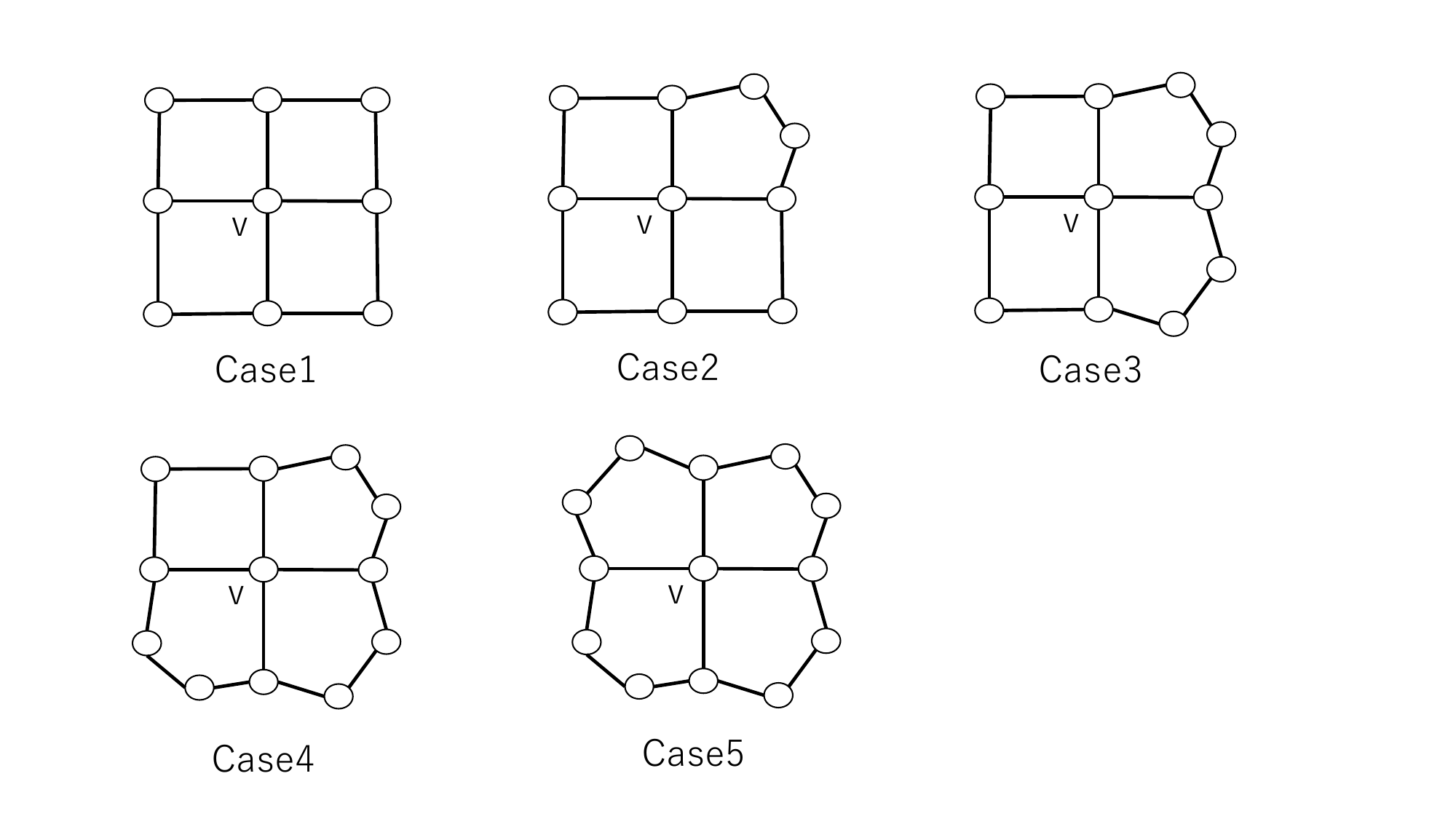}
\caption{Illustrations of a 4-vertex $v$ that is not incident to any 3-face.}
\label{fig_4-face}
\end{center}
\end{figure}
Case1: The 4-vertex $v$ is incident to four 4-faces.
By R9, the final charge is $2 - 4 \times \frac{1}{2} = 0$.
Case2: The 4-vertex $v$ is incident to three 4-faces and one 5-face.
By R9 and R10, the final charge is $2 - 3 \times \frac{1}{2} - \frac{1}{3} = \frac{1}{6}$.
Case3: The 4-vertex $v$ is incident to two 4-faces and two 5-faces.
By R9 and R10, the final charge is $2 - 2 \times \frac{1}{2} - 2 \times \frac{1}{3} = \frac{1}{3}$.
Case4: The 4-vertex $v$ is incident to one 4-face and three 5-faces.
By R9 and R10, the final charge is $2 - \frac{1}{2} - 3 \times \frac{1}{3} = \frac{1}{2}$.
Case5: The 4-vertex $v$ is incident to four 5-faces.
By R10, the final charge is $2 - 4 \times \frac{1}{3} = \frac{2}{3}$.
\end{itemize}
\item[5-vertex:] By Lemma~\ref{lem10}, a 5-vertex $v$ is $t(v) \leq 4$.
Therefore we divide the case by the value of $t(v)$.
\begin{itemize}
\item The case $t(v) = 4$ 

By Lemma~\ref{lem11}, the 5-vertex $v$ is incident to four (5,5,5)-faces.
By R1 and R11, the final charge is $4 - 4 \times 1 = 0$.
\item The case $t(v) = 3$

We further divide the case by the number of (5,5,5)-faces which is incident to $v$.
\begin{itemize}
\item The case where $v$ is incident to three (5,5,5)-faces.

\quad By Lemma~\ref{lem12}, the 5-vertex $v$ is not incident to two (5,5,5,3)-faces.
By Lemma~\ref{lem13}, the 5-vertex $v$ is not incident to one (5,5,5,3)-face and one $(4^+,4^+,4^+,4^+)$-face.
Thus, the worst situation is that the other two faces are one (5,5,5,3)-face and one 5-face or two $(4^+,4^+,4^+,4^+)$-faces.
By R1, R4 and R6, the final charge is at least $4 - 3 \times 1 - \frac{2}{3} - \frac{1}{3} = 0$.
By R1 and R5, the final charge is at least $4 - 3 \times 1 - 2 \times \frac{1}{2} = 0$.

\item The case where $v$ is incident to two (5,5,5)-faces and one (5,5,4)-face.

\quad By Lemma~\ref{lem15}, the 5-vertex $v$ is not incident to any (5,5,5,3)-faces.
By Lemma~\ref{lem16}, the 5-vertex $v$ is not incident to two $(4^+,4^+,4^+,4^+)$-faces.
Thus, the worst situation is that the other two faces are one 5-face and one $(4^+,4^+,4^+,4^+)$-face.
By R1, R2, R5 and R6, the final charge is at least $4 - 2 \times 1 - \frac{7}{6} - \frac{1}{2} - \frac{1}{3} = 0$.
\item The case where $v$ is incident to two (5,5,5)-faces and one (5,4,4)-face.

\quad By Lemma~\ref{lem14}, the 5-vertex $v$ is not incident to any 4-faces.
Thus, the worst situation is that the other two faces are two 5-faces.
By R1, R3 and R6, the final charge is at least $4 - 2 \times 1 - 1 - 2 \times \frac{1}{3} = \frac{1}{3}$.
\item The case where $v$ is incident to one (5,5,5)-face.

\quad By Lemma~\ref{lem17}, the 5-vertex $v$ is not incident to any 4-faces, and other two 5-faces are (5,5,4)-faces.
Thus, the worst situation is that the other two faces are two 5-faces.
By R1, R2 and R6, the final charge is at least $4 - 1 - 2 \times \frac{7}{6} - 2 \times \frac{1}{3} = 0$.
\item The case where $v$ is not incident to (5,5,5)-face.

\quad It is not possible for such a case. (see Figure~\ref{fig_a2}.)
\end{itemize}
\item The case $t(v) = 2$

We further divide the case by the number of (5,5,5)-faces which is incident to $v$.
\begin{itemize}
\item The case where $v$ is incident to two (5,5,5)-faces.

\quad The worst situation is that the other three faces are all (5,5,5,3)-faces.
By R1 and R4, the final charge is at least $4 - 2 \times 1 - 3 \times \frac{2}{3} = 0$.
\item The case where $v$ is incident to one (5,5,5)-face and one (5,5,4)-face.

\quad The worst situation is that the other three faces are two (5,5,5,3)-faces and one $(4^+,4^+,4^+,4^+)$-face.
By R1, R2, R4 and R5, the final charge is at least $4 - 1 - \frac{7}{6} - 2 \times \frac{2}{3} - \frac{1}{2} = 0$.
\item The case where $v$ is incident to one (5,5,5)-face and one (5,4,4)-face.

\quad By Lemma~\ref{lem18}, two of the other three faces are not 4-faces.
Thus, the worst situation is that the other three faces are one (5,5,5,3)-face and two 5-faces.
By R1, R3, R4 and R6, the final charge is at least $4 - 1 - 1 - \frac{2}{3} - 2 \times \frac{1}{3} = \frac{2}{3}$.
\item The case where $v$ is not incident to (5,5,5)-face.

\quad By Lemma~\ref{lem19}, if $v$ is incident to two (5,5,4)-faces, then two of the other three faces are not (5,5,5,3)-faces.
Thus, the worst situation is that the five faces are two (5,5,4)-faces, one (5,5,5,3)-face and two $(4^+,4^+,4^+,4^+)$-faces.
By R2, R4 and R5, the final charge is at least $4 - 2 \times \frac{7}{6} - \frac{2}{3} - 2 \times \frac{1}{2} = 0$.
\end{itemize}

\item The case $t(v) = 1$

We further divide the case by the face which is incident to $v$.
\begin{itemize}
\item The case where $v$ is incident to (5,5,5)-face.

\quad The worst situation is that the other four faces are all (5,5,5,3)-faces.
By R1 and R4, the final charge is at least $4 - 1 - 4 \times \frac{2}{3} = \frac{1}{3}$.
\item The case where $v$ is incident to (5,5,4)-face.

\quad The worst situation is that the other four faces are three (5,5,5,3)-faces and one $(4^+,4^+,4^+,4^+)$-face.
By R2, R4 and R5, the final charge is at least $4 - \frac{7}{6} - 3 \times \frac{2}{3} - \frac{1}{2} = \frac{1}{3}$.
\item The case where $v$ is incident to (5,4,4)-face.

\quad The worst situation is that the other four faces are two (5,5,5,3)-faces and two $(4^+,4^+,4^+,4^+)$-faces.
By R3, R4 and R5, the final charge is at least $4 - 1 - 2 \times \frac{2}{3} - 2 \times \frac{1}{2} = \frac{2}{3}$.
\end{itemize}
\item The case $t(v) = 0$

The worst situation is that $v$ is incident to five (5,5,5,3)-faces.
By R4, the final charge is at least $4 - 5 \times \frac{2}{3} = \frac{2}{3}$.
\end{itemize}
\end{enumerate}
Based on the above, the final charge of all vertices and faces are nonnegative, which is a contradiction.
Thus, Theorem~\ref{main} holds.

\bibliographystyle{elsarticle-num-names} 
\bibliography{ref}

\appendix
\section{Appendix}
Let $G$ be a minimum counterexample with minimum $|V(G)| + |E(G)|$ to Theorem \ref{main}.
Let $v$ is a 5-vertex with $N_G(v) = \{v_1,v_2,v_3,v_4,v_5\}$.
We show reducible configurations depending on the value of $t(v)$.
Red number in the figure indicate the degree of the vertex.
\subsection{$t(v) = 3$}
\begin{figure}[ht]
\begin{center}
\includegraphics[scale=0.25]{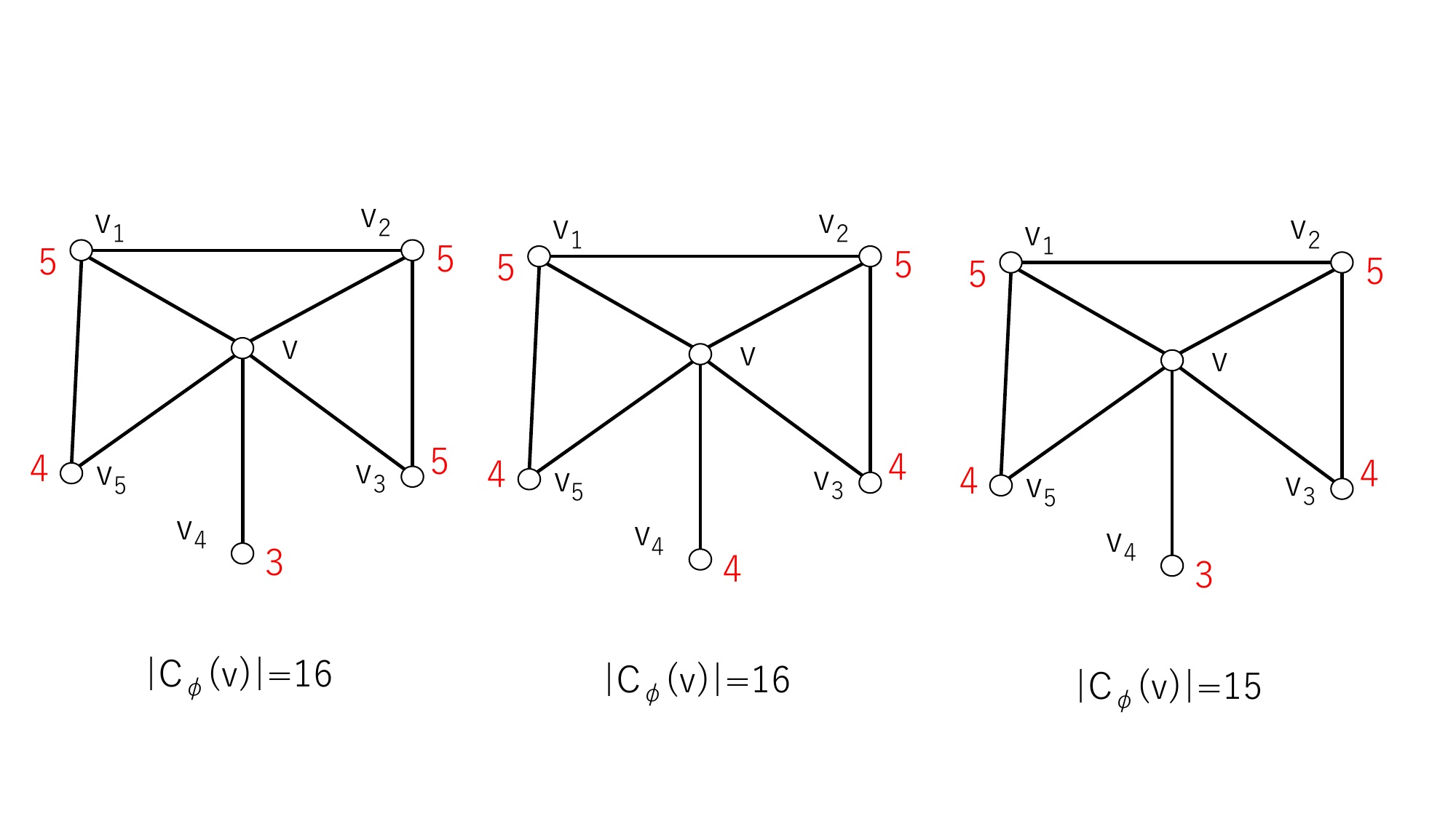}
\caption{Illustrations for case1 of $t(v) = 3$.}
\label{fig_a1}
\end{center}
\end{figure}

\begin{figure}[ht]
\begin{center}
\includegraphics[scale=0.25]{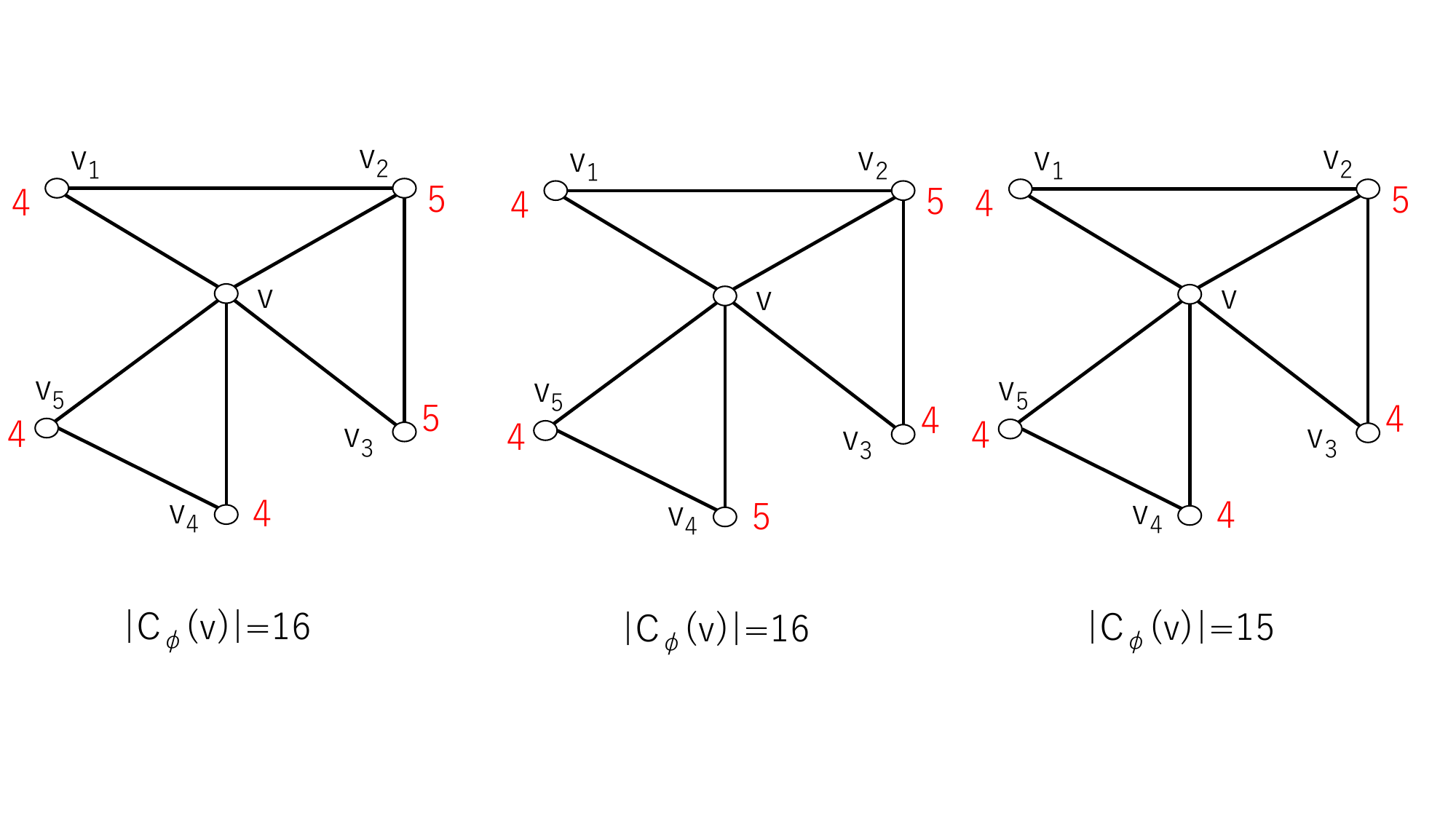}
\caption{Illustrations for case2 of $t(v) = 3$.}
\label{fig_a2}
\end{center}
\end{figure}
We have following two cases.

Case1: 3-face = $[vv_1v_2],[vv_2v_3],[vv_5v_1]$.
Case2: 3-face = $[vv_1v_2],[vv_2v_3],[vv_4v_5]$.
\begin{itemize}
\item Case1: (See Figure~\ref{fig_a1}.)
Let $G' = G - \{v\} + v_3v_4 + v_4v_5$.
By the minimality of $G$, $G'$ has a 2-distance 17-coloring $\phi'$.
Let $\phi$ be a coloring of $G$ such that every vertex in $V(G)$, except for $v$, is colored using $\phi'$.
All graphs satisfy $|C_\phi(v)| < 17$.
If $v$ is colored with $\phi(v) \in C \setminus C_\phi(v)$, then there exists a coloring $\phi$ of $G$ such that $\chi_2(G) \leq 17$, which is a contradiction.
\item Case2: (See Figure~\ref{fig_a2}.)
Let $G' = G - \{v\} + v_3v_4 + v_1v_5$.
By the minimality of $G$, $G'$ has a 2-distance 17-coloring $\phi'$.
Let $\phi$ be a coloring of $G$ such that every vertex in $V(G)$, except for $v$, is colored using $\phi'$.
All graphs satisfy $|C_\phi(v)| < 17$.
If $v$ is colored with $\phi(v) \in C \setminus C_\phi(v)$, then there exists a coloring $\phi$ of $G$ such that $\chi_2(G) \leq 17$, which is a contradiction.
\end{itemize}
\subsection{$t(v) = 2$}

\begin{figure}[ht]
\begin{center}
\includegraphics[scale=0.25]{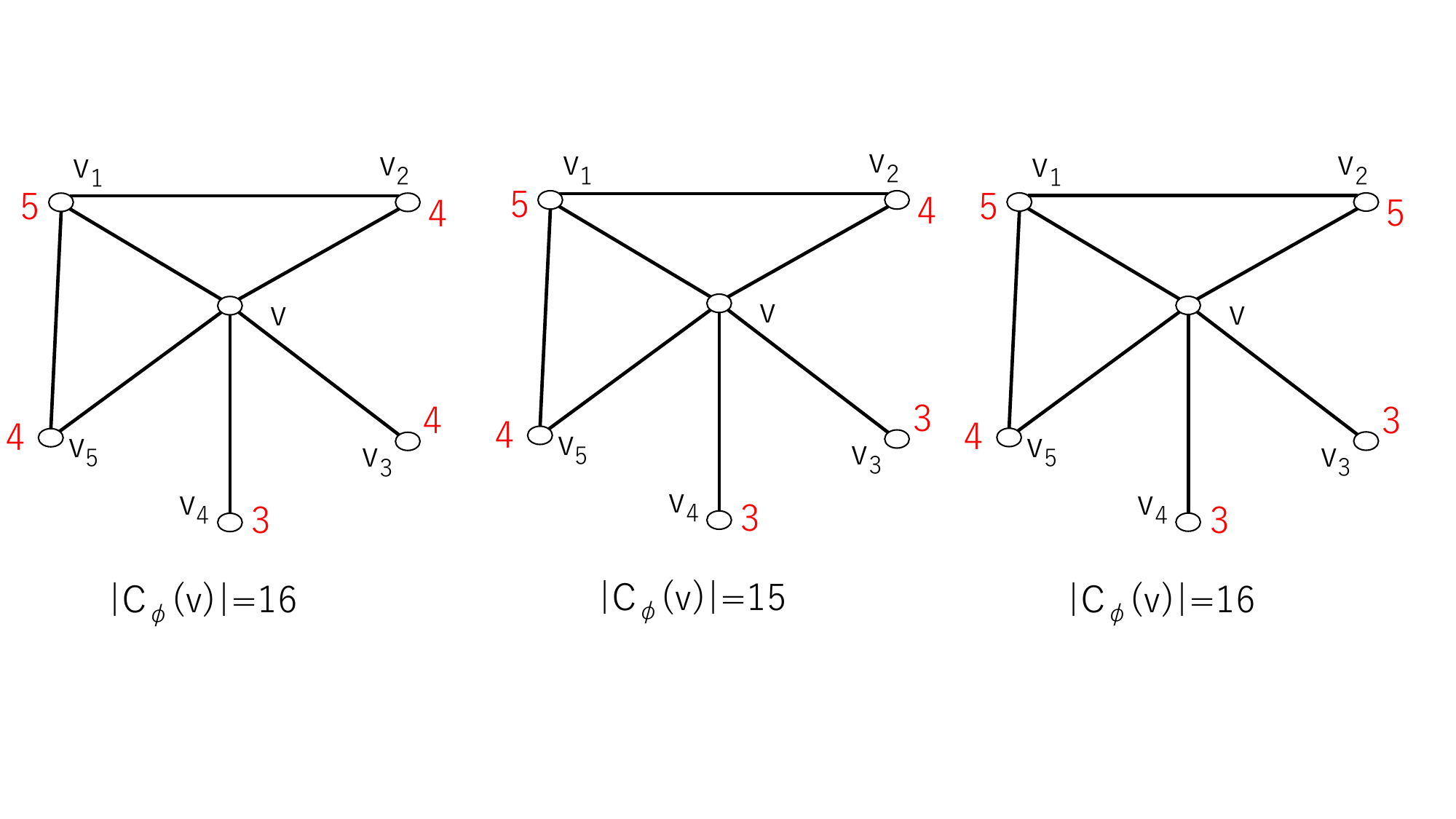}
\caption{Illustrations for case1 of $t(v) = 2$.}
\label{fig_a3}
\end{center}
\end{figure}

\begin{figure}[ht]
\begin{center}
\includegraphics[scale=0.25]{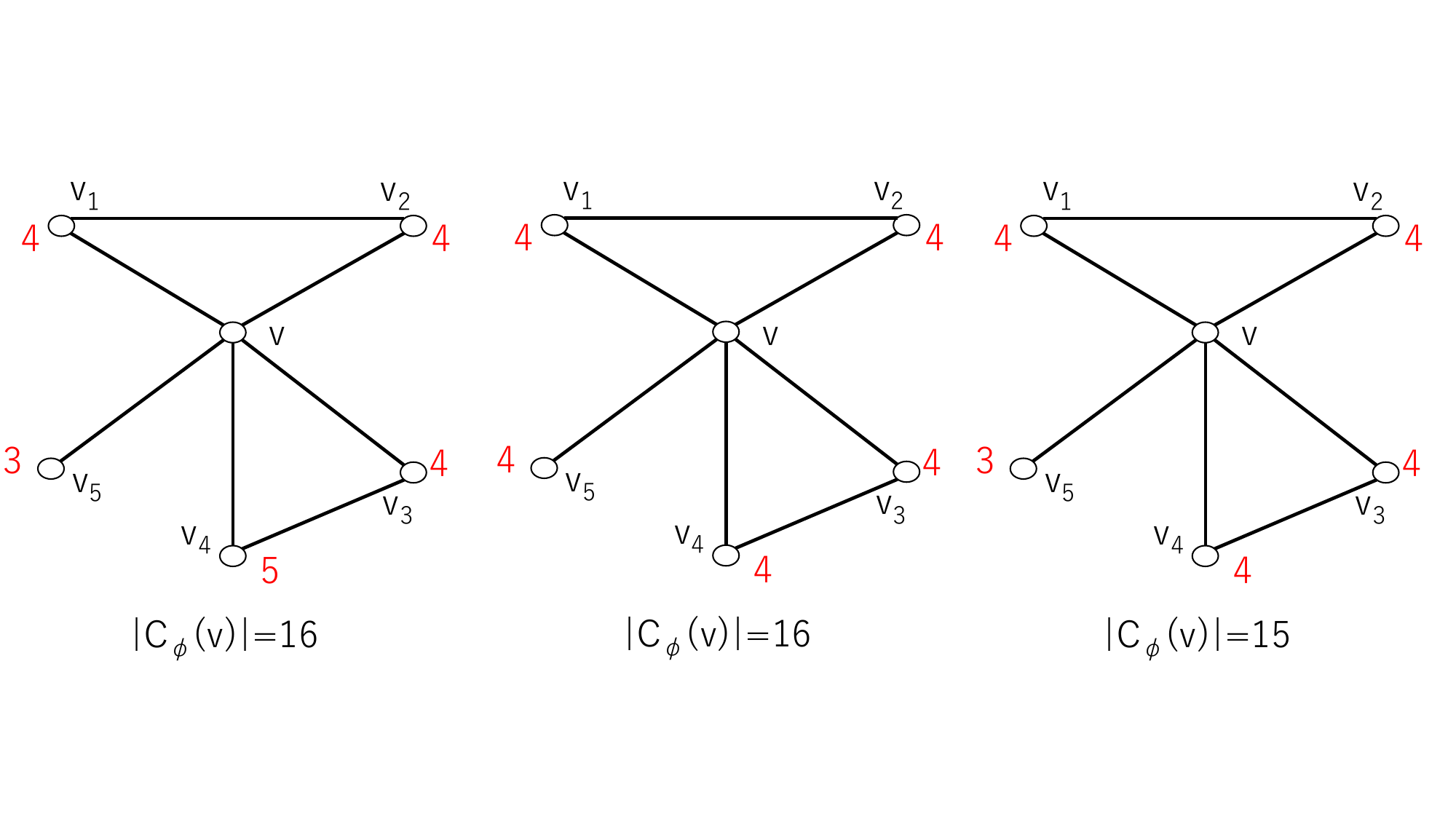}
\caption{Illustrations for case2 of $t(v) = 2$.}
\label{fig_a4}
\end{center}
\end{figure}

We have following two cases.
Case1: 3-face = $[vv_1v_2],[vv_5v_1]$.
Case2: 3-face = $[vv_1v_2],[vv_3v_4]$.
\begin{itemize}
\item Case1: (See Figure~\ref{fig_a3}.)
Let $G' = G - \{v\} +v_2v_3 + v_3v_4 + v_4v_5$.
By the minimality of $G$, $G'$ has a 2-distance 17-coloring $\phi'$.
Let $\phi$ be a coloring of $G$ such that every vertex in $V(G)$, except for $v$, is colored using $\phi'$.
All graphs satisfy $|C_\phi(v)| < 17$.
If $v$ is colored with $\phi(v) \in C \setminus C_\phi(v)$, then there exists a coloring $\phi$ of $G$ such that $\chi_2(G) \leq 17$, which is a contradiction.
\item Case2: (See Figure~\ref{fig_a4}.)
Let $G' = G - \{v\} + v_2v_3 + v_4v_5 + v_1v_5$.
By the minimality of $G$, $G'$ has a 2-distance 17-coloring $\phi'$.
Let $\phi$ be a coloring of $G$ such that every vertex in $V(G)$, except for $v$, is colored using $\phi'$.
All graphs satisfy $|C_\phi(v)| < 17$.
If $v$ is colored with $\phi(v) \in C \setminus C_\phi(v)$, then there exists a coloring $\phi$ of $G$ such that $\chi_2(G) \leq 17$, which is a contradiction.
\end{itemize}

\begin{figure}[ht]
\begin{center}
\includegraphics[scale=0.25]{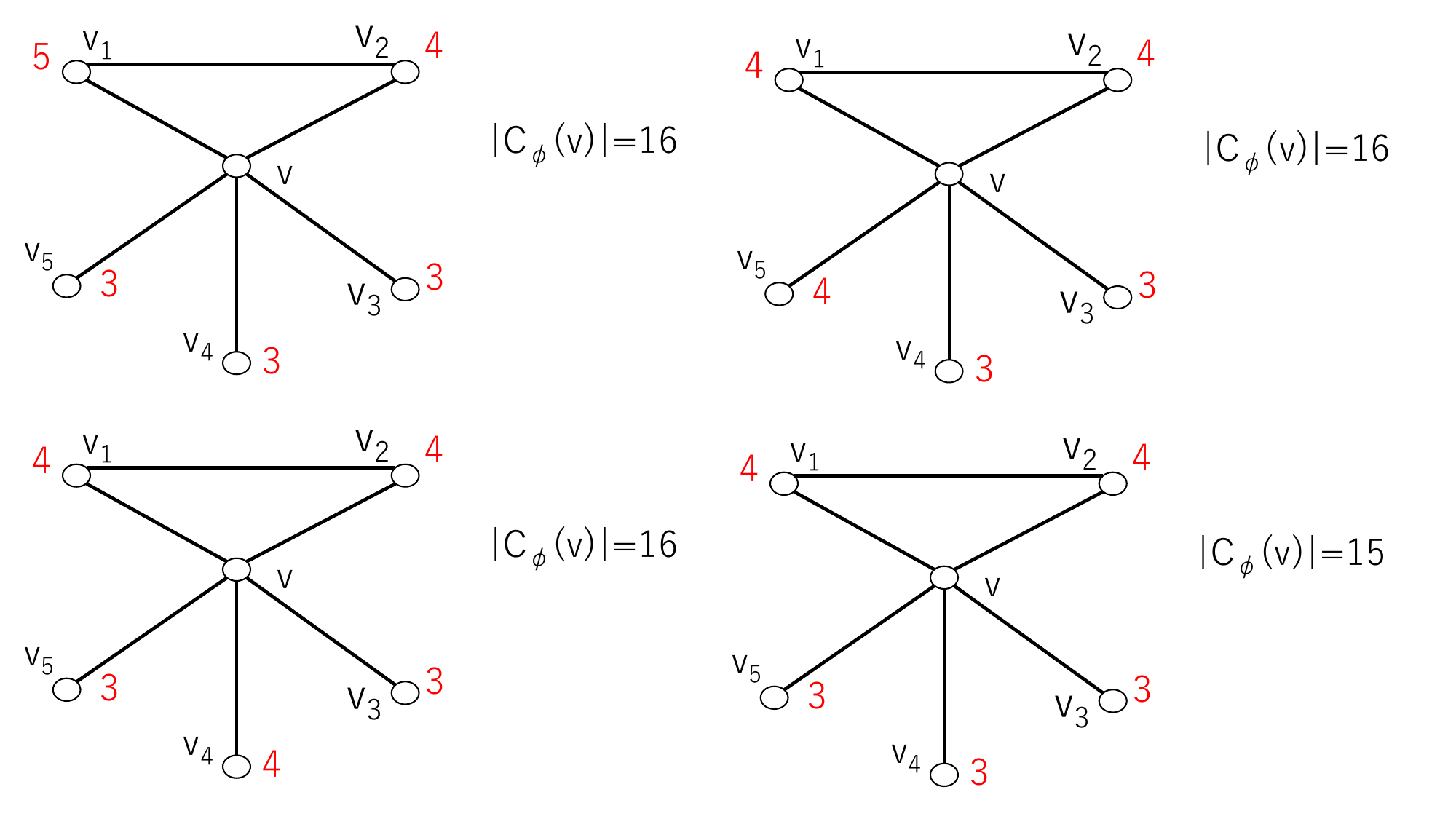}
\caption{Illustrations of $t(v) = 1$.}
\label{fig_a5}
\end{center}
\end{figure}

\subsection{$t(v) = 1$}
(See Figure~\ref{fig_a5}.)
Let $G' = G - \{v\} + v_2v_3 + v_3v_4 + v_4v_5 + v_1v_5$.
By the minimality of $G$, $G'$ has a 2-distance 17-coloring $\phi'$.
Let $\phi$ be a coloring of $G$ such that every vertex in $V(G)$, except for $v$, is colored using $\phi'$.
All graphs satisfy $|C_\phi(v)| \leq 17$.
If $v$ is colored with $\phi(v) \in C \setminus C_\phi(v)$, then there exists a coloring $\phi$ of $G$ such that $\chi_2(G) \leq 17$, which is a contradiction.
\end{document}